\documentclass[12pt]{amsart}
\usepackage{amsmath,amssymb, xypic,eepic,verbatim,amscd}
\numberwithin{equation}{section}

\newcommand\bull{\sssize{\bullet}}
\newcommand\home{\operatorname{Hom}}
\newcommand\spec{\operatorname{Spec}}
\newcommand\Spann{\operatorname{Span}}

\newcommand{\letof}[1]{\stackrel{#1}{\longrightarrow}}
\newcommand\mf{\mathcal{F}}
\newcommand\mk{\mathcal{K}}
\newcommand\ram{R}
\newcommand\me{\mathcal{E}}
\newcommand\mt{\mathcal{T}}

\newcommand\mrr{\mathcal{R}}
\newcommand\Fl{\operatorname{Fl}}
\newcommand\mv{\mathcal{V}}
\newcommand\tensor{\otimes}
\newcommand\ml{\mathcal{L}}
\newcommand\pl{P}
\newcommand\tpl{\mathcal{P}}
\newcommand\ma{\mathcal{A}}
\newcommand\shom{\mathcal{H}{\it{om}}}
\newcommand\codim{\operatorname{codim}}
\newcommand\mb{\mathcal{B}}
\newcommand\mpp{\mathcal{P}}
\newcommand\GL{\operatorname{GL}}
\newcommand\mj{\mathcal{J}}
\newcommand\Fone{J}
\newcommand\Ftwo{\mj}
\newcommand\mg{\mathcal{G}}
\newcommand\mh{\mathcal{H}}

\newcommand\im{\operatorname{im}}
\newcommand\mw{\mathcal{W}}
\newcommand\mc{\mathcal{C}}
\newcommand\mq{\mathcal{Q}}
\newcommand\rk{\operatorname{rk}}
\newcommand\ms{\mathcal{S}}

\newcommand\mi{\mathcal{I}}

\newcommand\Gr{\operatorname{Gr}}

\newcommand\muu{\mathfrak{U}}
\newcommand\mukv{\mathfrak{U}_{\mk}(V)}
\newcommand\muiw{\mathfrak{U}_{\mi}(W)}
\newcommand\mhh{\mathfrak{H}}
\newcommand{\leto}[1]{\stackrel{#1}{\to}}
\newcommand\hoo{\home_{\mi}(V,Q,\mf,\mg)}

\newcommand\ham{\mhh_{\mi}(V,Q,\mk)}
\newcommand\capn{P}
\newcommand\smalln{p}

\newcommand\hmk{H}

\newtheorem{theorem}{Theorem}[section]

\newtheorem{corollary}[theorem]{Corollary}

\newtheorem{proposition}[theorem]{Proposition}

\newtheorem{lemma}[theorem]{Lemma}
\newtheorem{clm}[theorem]{Claim}
\newtheorem{definition}[theorem]{Definition}

\newtheorem{remark}[theorem]{{ Remark}}
\newtheorem{notation}[theorem]{Notation}
\theoremstyle{remark}

\newtheorem{example}[theorem]{\bf Example}

\begin{document}
\title[Horn and Saturation conjectures]{Geometric Proofs of Horn and Saturation conjectures}
\author{Prakash Belkale}
\address{Department of Mathematics\\ UNC-Chapel Hill\\ CB \#3250, Phillips Hall
\\ Chapel Hill,  NC 27599}
\email{belkale@email.unc.edu}

\begin{abstract}
We provide a geometric proof of the Schubert calculus interpretation  of the Horn conjecture, and show how the saturation conjecture  follows from it. The geometric proof gives a  strengthening of Horn and saturation conjectures. We also establish transversality theorems for Schubert calculus in non-zero characteristic.
\end{abstract}
\footnote{The author was partially supported by NSF grant DMS-0300356.}

\maketitle

\section*[intro]{Introduction}
One of the consequences of the work of A. Klyachko ~\cite{Klyachko}, and the saturation theorem of A. Knutson and T. Tao ~\cite{KT} is that one can tell when a product of Schubert classes in a given Grassmannian is nonzero by writing down a series of inequalities coming from knowing the answer to the same question for smaller Grassmannians. The proofs involve the connections of this problem with representation theory, combinatorics and geometric invariant theory. In ~\cite{bourful}, W. Fulton proposed the challenge of finding a geometric proof of this property. We will refer to this property as Geometric Horn (GH) because of its close relations to
  a 1962 conjecture of A. Horn ~\cite{Horn} on the characterization of the possible eigenvalues of a sum of Hermitian matrices in terms of the eigenvalues of the summands. The original conjecture of Horn was proved by the combined works of Klyachko ~\cite{Klyachko} and Knutson-Tao ~\cite{KT}. We refer the reader to Fulton's  article ~\cite{ful} for a discussion of these problems.

In this paper we give a geometric proof of a generalization of GH and obtain a
stronger form of the saturation problem as a consequence.

Following a suggestion of F. Sottile, we use the geometric proof of GH to  establish  transversality statements for Schubert calculus of Grassmannians in any characteristic (see Corollary ~\ref{sottilee}). A proof of transversality also appears in a recent paper of R. Vakil ~\cite{vakil}.

Let us first introduce some notation, and then recall the statement of GH, phrasing it in a slightly different form than in \cite{ful}.

Let $I$ be a subset of $[n]=\{1,\dots,n\}$ of cardinality $r$. Make the convention that a set $I$ of this form is always written in the form $I=\{i_1<\dots<i_r\}$ (if we have a collection of such subsets $I^j,$ $ j=1,\dots, s$ each of cardinality $r$ we will assume
that they are written as $I^j=\{i^j_1<\dots<i^j_r\}$ for $j=1,\dots,s$). Let $W$ be an $n$-dimensional vector space and
$$E_{\sssize{\bullet}}:\text{ }\{0\}=E_0\subsetneq E_1\subsetneq\dots\subsetneq E_n=W$$
 be  a complete flag of subspaces of $W$. Define $\Omega_{I}(E_{\sssize{\bullet}})\subseteq \Gr(r,W)$ to be
$$ \{ V \in \Gr(r,W) \mid \rk(V \cap
E_{i_{a}}) \geq a\text{ for }1\leq a \leq r\}.$$ We denote the class in the integral cohomology ring (or the Chow ring, see Section ~\ref{notate}) of
this subvariety by $\omega_{I}$. The codimension of $\Omega_{I}(E_{\bull})$ in $\Gr(r,n)$ is
$$\codim(\omega_I)=\sum_{a=1}^{r} (n-r +a -i_a).$$
We give a geometric proof of the following theorem (due to Klyachko, Knutson and Tao):

{\bf Geometric Horn}: Let $I^1,\dots,I^s$ be subsets of
$[n]=\{1,\dots,n\}$ each of cardinality $r$. Assume that
$$\sum_{j=1}^{s}\sum_{a=1}^{r} (n-r +a -i^j_a)= r(n-r).$$
The following are equivalent:
\begin{enumerate}
\item
$\prod_{j=1}^s\omega_{I^j}= \text{non-zero multiple of a class of a point in } H^{*}(\Gr(r,n)).$
\item For any $1 \leq d < r$ and any choice of subsets $K^{1},\dots,K^{s}$ each of cardinality $d$ of  $[r]=\{1,\dots,r\}$ such that
$\prod_{j=1}^s\omega_{K^j}=\text{non-zero multiple of the class of a point in }$ $H^{*}(\Gr(d,r)),$
 the  inequality
$\sum_{j=1}^s \sum_{a\in K^{j}}(n-r+a-i^j_a) \leq d(n-r)$
is valid.
\end{enumerate}
 Our geometric approach gives a stronger form of GH in which inductive conditions for a cup product to be non-zero are given:
\begin{theorem}\label{MainTe} Let $\mi=(I^1,\dots,I^s)$ be a $s$-tuple of subsets of
$[n]=\{1,\dots,n\}$ each of cardinality $r$.  The following are equivalent:
\begin{enumerate}
\item[(A)] $\prod_{j=1}^s\omega_{I^j}\neq 0  \text{ in } H^{*}(\Gr(r,n)).$
\item[(B)]  For any $0 < d \leq  r$ and any $s$-tuple $\mk=(K^{1},\dots,K^{s})$ of subsets of $[r]=\{1,\dots,r\}$ each of cardinality $d$ such that $\prod_{j=1}^s\omega_{K^j}\neq 0\text{ in }  H^{*}(\Gr(d,r))$, the inequality
\begin{equation}\nonumber\tag{$\dagger_{\mk}^{\mi}$}
   \sum_{j=1}^s \sum_{a\in K^{j}}(n-r+a-i^j_a)-d(n-r) \leq 0
\end{equation}
is valid.
\item[(C)]  For any $0 <d \leq r$ and any $s$-tuple  $\mk=(K^{1},\dots,K^{s})$ of subsets of cardinality $d$ of $[r]=\{1,\dots,r\}$ such that $\prod_{j=1}^s\omega_{K^j}=\text{class of a point in }  H^{*}(\Gr(d,r))$, the Inequality $(\dagger_{\mk}^{\mi})$ is valid.
\end{enumerate}
\end{theorem}
Notice that in Theorem ~\ref{MainTe}, there is a unique choice of $\mk$ in (B) (or in (C)) subject to the constraint $d=r$.  The corresponding inequality
$(\dagger_{\mk}^{\mi})$ is the codimension condition:
\begin{equation}\label{conditione}
\sum_{j=1}^{s}\sum_{a=1}^{r} (n-r +a -i^j_a)\leq r(n-r).
\end{equation}
 Because of the isomorphism between $\Gr(r,n)$ and
$\Gr(n-r,n)$, there is a  dual set of inequalities which characterize when a cup product of cycle classes of Schubert varieties in $\Gr(r,n)$ is nonzero. This is not the same set of inequalities, and it remains somewhat of a mystery why these give the same conditions.

Theorem ~\ref{MainTe} implies the saturation theorem of Knutson-Tao ~\cite{KT} in a stronger form. Recall that irreducible {\em polynomial} representations of $\GL(r)$ are indexed by sequences $\lambda=(\lambda_1\geq\dots\geq \lambda_r\geq 0)\in \Bbb{Z}^r$. Denote the representation corresponding to $\lambda$ by $V_{\lambda}$. Define Littlewood-Richardson  coefficients $c^{\lambda}_{\mu,\nu}$ by: $V_{\mu}\otimes V_{\nu}=\sum c^{\lambda}_{\mu,\nu}V_{\lambda}.$
The saturation theorem of Knutson-Tao states that for any positive integer $N$,
$$c^{\lambda}_{\mu,\nu}\neq 0\text{  if and only if  } c^{N\lambda}_{N\mu,N\nu}\neq 0.$$
Our work gives a generalization of the saturation theorem by  bounding the widths of the sequences involved (the extension of the saturation theorem to more than two factors was known previously):
\begin{theorem}\label{S1} Let $V_{\lambda}$, $V_{\mu}$, $\dots, V_{\nu}$ be irreducible representations of $\GL(r)$, $N$ and $\ell$  positive integers. The following are equivalent:
\begin{enumerate}
\item
There exist a $V_{\delta}$  with $\delta_1 \leq \ell$ and an inclusion of $\GL(r)$ representations
$V_{\delta}\subseteq  V_{\lambda}\otimes V_{\mu}\otimes\dots\otimes V_{\nu}$.
\item
There exist a $V_{\delta'}$  with $\delta'_1 \leq N\ell$ and an inclusion of $\GL(r)$ representations
$V_{\delta'}\subseteq  V_{N\lambda}\otimes V_{N\mu}\otimes\dots\otimes V_{N\nu}$.
\end{enumerate}
\end{theorem}

We indicate briefly the main idea of our proof: Assume that the ground field has characteristic zero. If general Schubert varieties intersect at a point, then by Kleiman transversality they intersect transversally there. Conversely one can detect if general Schubert varieties intersect by a tangent space calculation. The tangent space of a Grassmannian $\Gr(r,W)$ at a point $V$ is canonically isomorphic to $\home(V,W/V)$. We use the action of $GL(V)\times \GL(W/V)$ on $\home(V,W/V)$ (this action has only finitely many orbits) to study the intersection of the tangent spaces of general Schubert varieties passing through the point $V$.

 We use techniques from the  theory of parabolic bundles ~\cite{MS}. The intersection of tangent spaces of Schubert varieties at a point of intersection (in a Grassmannian) is analogous to the vector space of morphisms between two parabolic bundles. 

The following general position idea is used in a crucial way: In certain cases, induced structures in an intersection theory situation can be assumed to be ``generic''.  For example, if three general Schubert varieties in a Grassmannian $\Gr(r,W)$ meet at a point $V$, then $V$ gets three induced complete flags. Can we assume that these three induced flags on $V$ are ``generic enough'' for some other intersection theory calculation? We develop techniques to answer this kind of question by connecting it to  the question of irreducibility of certain parameter spaces (see Proposition ~\ref{neww}). 

In ~\cite{bell}, we use similar methods to prove a quantum (multiplicative) generalization of Horn's conjecture.

I would like to acknowledge with gratitude, generous help from W. Fulton in writing this paper. In addition, I would like to thank  P. Brosnan, S. Kumar,  J. Millson, M.V. Nori and F. Sottile for useful discussions.
 \subsection{Notation:}\label{notate}
 We make the following conventions:
\begin{itemize}
\item A  variety is a reduced and irreducible scheme over a  fixed algebraically closed field $\kappa$ of arbitrary characteristic. By a vector space, we mean a vector space over the field $\kappa$. The Chow ring of the Grassmannian $\Gr(r,n)$ will be denoted by $H^*(\Gr(r,n),\Bbb{Z})$ or simply $H^*(\Gr(r,n))$. This should not cause any confusion since the Chow ring of a Grassmannian does not depend on the characteristic and is isomorphic (via the cycle class map) to singular cohomology if the base field is $\Bbb{C}$.
\item The cycle class (in cohomology or in the Chow group) of a closed subvariety $Z$ of a smooth variety $X$ is denoted by $[Z]$. Recall that the class in cohomology is the Poincar\'{e}  dual of the homology class determined by $Z$.
\item We fix an integer $s\geq 1$.
\item For a vector space $W$, let $\Fl(W)$ denote the variety of complete flags on it. If $\me\in \Fl(W)^s$, we will assume that $\me$ is
written in the form $(E^1_{\bull},\dots,E^s_{\bull})$.
\item We use the notation $[n]=\{1,\dots,n\}$.
\item Let {{\bf P}} be a property that makes sense for closed points of a given variety $X$. We say that a generic point of $X$ satisfies {\bf{P}} (or, for generic $x\in X$, {\bf {P}} holds) if there is a nonempty Zariski open subset $U$ of $X$ such that {\bf {P}} holds for every $x\in U$.
\end{itemize}
\section{Preliminaries from Schubert calculus}\label{sec1}
Let $V$ be a point of intersection of Schubert varieties in a Grassmannian $\Gr(r,W)$. $V$ acquires additional structure, namely induced flags and this additional structure can be used consider Schubert calculus in the Grassmannians $\Gr(d,V)$. There is also a natural inclusion $\Gr(d,V)\hookrightarrow \Gr(d,W)$. We therefore see hints of a recursive structure in Schubert calculus. This section lays groundwork for such arguments and  proves (A)$\Rightarrow$ (B) in Theorem ~\ref{MainTe}. We begin by recalling basic facts from dimension theory.
\subsection{Basic dimension  and intersection theory}\label{faizfaiz}
 Let $X$ be a variety and  $Y,$ $Z$  locally closed subvarieties of $X$. Assume that $Z$ is locally defined by the vanishing of $p$ functions. By Krull's principal ideal theorem, each irreducible component of $Y\cap Z$ is of dimension $\geq \dim(Y)-p$.
  If $X$ is smooth, the diagonal  $\Delta_X\subseteq  X\times X$ is given locally by the vanishing of $\dim(X)$ functions, hence each irreducible  component of
$Y\cap Z=(Y\times Z)\cap \Delta_X$ is of dimension at least $\dim(Y)+\dim(Z)-\dim(X). $

Let $X$ be a smooth variety and $X_1,\dots,X_s$  locally closed subvarieties  of $X$. By considering the intersection of the main diagonal $\Delta_X$ with $X_1\times\dots\times X_s$ in $X^s$, one finds that each irreducible  component of $\cap_{j=1}^s X_j$ is of dimension at least
$\dim(X)-\sum_{j=1}^s \codim(X_j,X).$ An irreducible component $Z$ of the intersection of $X_1,$ $\dots,$ $X_s$ is said to be {\bf proper} if it is of dimension $\dim(X)-\sum_{j=1}^s \codim(X_j,X).$

 Suppose, in addition, that the subvarieties $X_1,\dots,X_s$ are smooth, and $x\in \cap_{j=1}^s X_j$. We say that $X_1,\dots,X_s$ meet {\bf transversally} at $x$ if the intersection of the tangent spaces $T(X_j)_x$ of $X_j$ at $x$ for $j=1,\dots,s$ (these are subspaces of the tangent space to $X$ at $x$) is of dimension
$$\rk(\cap_{j=1}^sT(X_j)_x)=\dim(X)-\sum_{j=1}^s\codim(X_j,X).$$
If $x$ is a transverse point of intersection of the varieties $X_1,$ $\dots,$ $X_s$, then there is a unique irreducible component $Z$ of the scheme-theoretic intersection $\cap_{j=1}^s X_j$ containing $x$. Furthermore,  $Z$ is a proper irreducible component of the intersection $\cap_{j=1}^s X_j$ and is smooth at $x$.

We recall the  following result of Kleiman ~\cite{kl} on proper and transverse intersections:
\begin{proposition}\label{kleiman}
Suppose that a connected algebraic group $G$ acts transitively on a smooth variety $X$. Let $X_1$, $\dots,X_s$ be locally closed subvarieties of $X$. Then, there exists a non empty open subset $U\subseteq  G^s$ such that for $(g_1,\dots,g_s)\in U$,
\begin{enumerate}
\item Each irreducible component of the intersection of $g_1 X_1$, $\dots$, $g_s X_s$ is proper.
\item $\bigcap_{j=1}^s g_j X_j$ is dense in $\bigcap_{j=1}^s g_j \bar{X}_j$ (which could be empty).
\end{enumerate}
If the base field is of characteristic zero and $X_1$, $\dots$, $X_s$ are smooth varieties, we can find such a
$U$ with the  additional property that for $(g_1,\dots,g_s)\in U$, $g_1X_1,$ $\dots,$ $g_s X_s$ meet transversally at each point of their intersection.
\end{proposition}
\begin{proof}
We include a proof of the density statement (the rest of the conclusion is standard, see ~\cite{kl}).

Let  $Y_j=\bar{X}_j\smallsetminus X_j$ for $j=1,\dots s$. Let $U$ be a nonempty open subset of $G^s$  such that for $(g_1,\dots,g_s)\in U$, the following intersections are proper:
\begin{enumerate}
\item $\cap_{j=1}^s g_j \bar{X}_j$.
\item For $\ell\in \{1,\dots,s\}$, $\{\cap_{j\in\{1,\dots, s\}\smallsetminus \{\ell\}} g_j\bar{X}_j\}\cap g_{\ell} Y_{\ell}.$
\end{enumerate}
 For $(g_1,\dots,g_s)\in U$ and $\ell \in \{1,\dots,s\}$, each irreducible component of the intersection
$\{\cap_{j\in\{1,\dots, s\}\smallsetminus \{\ell\}} g_j\bar{X}_j\}\cap g_{\ell} Y_{\ell}$
is therefore of dimension strictly less than that of  each irreducible component of $\bigcap_{j=1}^s g_j \bar{X}_j$ ($\dim(Y_{\ell})<\dim(X_{\ell}),$ $ \ell=1,\dots,s$). This proves the density statement.
\end{proof}
The following is a standard result from intersection theory:
\begin{proposition}\label{inttheorie}
Let $X$ be a smooth projective variety and  $X_1$, $\dots,X_s$ closed subvarieties of $X$. Suppose that the cup product of the cycle classes $\prod_{j=1}^s[X_j]\neq 0\in H^*(X)$ (or the Chow ring). Then, $\cap_{j=1}^s X_j\neq \emptyset$.
\end{proposition}
Proposition ~\ref{inttheorie} is best understood in terms of Fulton and MacPherson's construction (see ~\cite{int}, Chapter 8) of a rational equivalence class of algebraic cycles  on $\cap_{j=1}^s X_j$ which represents the intersection product $X_1\cdot \ldots \cdot X_s$ (in the Chow ring of $X$). Hence if the product of cycle classes is nonzero, the intersection $\cap_{j=1}^s X_j$ is nonempty.
\subsection{Schubert cells in Grassmannians}
Let $I\subseteq [n]$ be a subset of cardinality $r$. Let $E_{\sssize{\bullet}}$ be a complete flag in an $n$-dimensional vector space $W$. Define the Schubert cell $\Omega^o_{I}(E_{\sssize{\bullet}})\subseteq \Gr(r,W)$ by
$$\Omega^o_I(E_{\sssize{\bullet}})=\{V\in \Gr(r,W)\mid \rk(V\cap E_u)=a \text{ for } i_a\leq u< i_{a+1},\text{ }a=0,\dots,r\}$$
where $i_0$ is defined to be $0$ and $i_{r+1}=n$.
$\Omega^o_{I}(E_{\sssize{\bullet}})$ is smooth and is an open dense subset of ${\Omega}_{I}(E_{\sssize{\bullet}})$. For a fixed complete flag on $W$, it is easy to see that (~\cite{ful2}, \S 1) every $r$-dimensional vector subspace  belongs to a unique Schubert cell.
\subsection{Induced flags}
Suppose that $W$ is an $n$-dimensional vector space and $V\subseteq  W$ an $r$-dimensional subspace. Let $E_{\bull}$ be a complete flag on $W$. This induces a complete flag on $V$ and a complete flag on $W/ V$ by intersecting $E_{\bull}$ with $V$ and by projection $p:W\to W/ V$  respectively. We denote these by $E_{\bull}(V)$ and $E_{\bull}(W/ V)$ respectively. Explicitly, if $V\in \Omega^o_I(E_{\bull})$ and $[n]\smallsetminus I=\{\alpha(1)<\dots<\alpha(n-r)\}$, then $E_a(V)=E_{i_a}\cap V,\ a=1,\dots,r$
and $E_b(W/ V)=p(E_{\alpha(b)})\text{, }b=1,\dots,n-r.$
Given an ordered collection of flags $\me\in \Fl(W)^s$ we obtain
ordered collections of flags $\me(V)\in \Fl(V)^s$ and $\me(W/ V)\in \Fl(W/ V)^s$
by performing the above operations in each coordinate factor.

The following lemma follows from a direct calculation (see ~\cite{ful2}, Lemma 2 (i)).
\begin{lemma}\label{falcon}
Let $W$ be an $n$-dimensional vector space. Suppose $E_{\bull}\in \Fl(W)$ and  $S\subseteq  V\subseteq  W $ are subspaces with $\rk(V)=r$ and $\rk(S)=d$. Let $I$ be the unique subset of $[n]$ of cardinality $r$  such that $V\in \Omega^o_{I}(E_{\bull})\subseteq \Gr(r,W)$, and $K$ the unique subset of $[r]$ of cardinality $d$  such that $S\in \Omega^o_{K}(E_{\bull}(V))\subseteq \Gr(d,V)$. Set $L=\{i_a\mid a\in K\}$. Then,
 $S\in\Omega^o_L(E_{\bull})\subseteq \Gr(d,W).$
\end{lemma}
\subsection{Proof of (A)$\Rightarrow$(B) in Theorem ~\ref{MainTe}}\label{First}
\begin{definition}\label{defone}Let $V$ be a $r$-dimensional subspace of an
$n$-dimensional vector space $W$, and $\me\in \Fl(W)^s$. Let $I^1,\dots, I^s$ be the unique subsets of $[n]$ each
of cardinality $r$ such that $V\in \Omega^o_{I^j}(E^j_{\bull})$ for $j=1,\dots,s$. Define $\dim(V,W,\me)$ to be the expected dimension of the intersection $\cap_{j=1}^s\Omega^o_{I^j}(E^j_{\bull})$. That is,
$$\dim(V,W,\me)=\dim(\Gr(r,n))-\sum_{j=1}^s \codim(\omega_{I^j})$$
$$= r(n-r)-\sum_{j=1}^s\sum_{a=1}^r (n-r+a-i^j_a).$$
\end{definition}\pagebreak
 \begin{proposition}\label{l2}
Let $W$ be an $n$-dimensional vector space.
\begin{enumerate}
\item Suppose $\me\in \Fl(W)^s$ and  $S\subseteq  V\subseteq  W $ are subspaces with $\rk(V)=r$ and $\rk(S)=d$. Let $\mi=(I^1,\dots,I^s)$ be the unique $s$-tuple of subsets of $[n]$ each of cardinality $r$  such that $V\in \cap_{j=1}^s\Omega^o_{I^j}(E^j_{\bull})$, and $\mk=(K^1,\dots,K^s)$  the unique $s$-tuple of subsets of $[r]$ each of cardinality $d$  such that $S\in \cap_{j=1}^s\Omega^o_{K^j}(E^j_{\bull}(V))$.  Then,
\begin{equation}\label{E100}
\dim(S,V,\me(V))-\dim(S,W,\me)=\sum_{j=1}^s\sum_{a\in K^j}(n-r+a-i^j_a)-d(n-r).
\end{equation}
\item There exists a nonempty  open subset $U$ of $\Fl(W)^s$ depending only on $W$ such that if $\me$, $S$, $V$, $\mi$ and $\mk$ are as in (1) with the additional condition that $\me\in U$, then
 \begin{equation}\label{E101}
\dim(S,V,\me(V))-\dim(S,W,\me)\leq 0
\end{equation}
 and  Inequality $(\dagger_{\mk}^{\mi})$ (which makes sense even if $d=0$, and in that case is just the inequality $0\leq 0$) holds.
\end{enumerate}
\end{proposition}
\begin{proof}
For (1), define $L^j=\{i^j_a\mid a\in K^j\}$ for $j=1\dots s$ which is a subset of $[n]$ of cardinality $d$. By Lemma ~\ref{falcon},
\begin{equation}\label{april}
S\in \bigcap_{j=1}^s \Omega^o_{K^j}(E^j_{\bull}(V))\subseteq  \bigcap_{j=1}^s \Omega^o_{L^j}(E^j_{\bull})\subseteq  \Gr(d,n).
\end{equation}
Hence, $$\dim(S,V,\me(V))=d(r-d)-\sum_{j=1}^s\sum_{t=1}^d(r-d+t-k^j_t)$$ and
$$\dim(S,W,\me)=d(n-d)-\sum_{j=1}^s\sum_{t=1}^d(n-d+t-i^j_{k^j_t}).$$
and we obtain Equation ~\ref{E100} by subtraction.

For (2), note that each irreducible component of $\cap_{j=1}^s \Omega^o_{K^j}(E^j_{\bull}(V))$ is of dimension at least
$\dim(S,V,\me(V))=\dim(Gr(d,r))-\sum \codim(\omega_{K^j})$. For generic $\me$, each irreducible
component of $\cap_{j=1}^s \Omega^o_{L^j}(E^j_{\bull})$ is of dimension exactly $\dim(S,W,\me)=\dim(Gr(d,n))-\sum
\codim(\omega_{L^j}).$

Inclusion ~\ref{april} therefore gives Inequality ~\ref{E101}. Note that Equation ~\ref{E100} and Inequality ~\ref{E101} imply Inequality $(\dagger_{\mk}^{\mi})$.
\end{proof}
Proposition ~\ref{l2} implies that $(A)\text{}\Rightarrow(B)$ in Theorem ~\ref{MainTe}. To see this assume that  $\prod_{j=1}^s\omega_{I^j}\neq 0\in H^*(\Gr(r,n))$, and $\prod_{j=1}^s\omega_{K^j}\neq 0\in H^*(Gr(d,r))$ as in the theorem. Let $W$ be an $n$ dimensional vector space and $\me\in \Fl(W)^s$ a generic point.

Pick $V\in\cap_{j=1}^s \Omega^o_{I^j}(E^j_{\bull})$. Since, $\prod_{j=1}^s \omega_{K^j}\neq 0 $, by Proposition ~\ref{inttheorie}, $\cap_{j=1}^s \Omega_{K^j}(E^j_{\bull}(V))\neq \emptyset$. Pick  $S\in \cap_{j=1}^s\Omega_{K^j}(E^j_{\bull}(V))$. Let $\mt=(T^1,\dots,T^s)$ be the unique $s$-tuple of subsets of $[r]$ such that  $S\in \cap_{j=1}^s\Omega^o_{T^j}(E^j_{\bull}(V))$. For $j=1,\dots,s$, $S\in \Omega_{K^j}(E^j_{\bull}(V))$ and $S\in \Omega^o_{T^j}(E^j_{\bull}(V))$, and therefore for $b=1,\dots,d$,
\begin{equation}\label{jjj}
t^j_b\leq k^j_b.
\end{equation}
 By Proposition ~\ref{l2}, (2) we obtain Inequality $(\dagger_{\mt}^{\mi})$:
$$\sum_{j=1}^s\sum_{a\in T^j}(n-r+a-i^j_a)-d(n-r)\leq 0.$$
The numbers $n-r+a-i^j_a$ weakly decrease as $a$ increases (with $j$ fixed). Therefore, from Inequality ~\ref{jjj} it follows that for $j=1,\dots,s$ and $b=1,\dots,d$,
$$(n-r+t^j_b-i^j_{t^j_b})\geq (n-r+k^j_b-i^j_{k^j_b}).$$
Inequality $(\dagger_{\mk}^{\mi})$ therefore follows from Inequality $(\dagger_{\mt}^{\mi})$.

The author learned from Sottile that the above approach towards the implication $(A)\Rightarrow(B)$ in Theorem ~\ref{MainTe} was known to him. The above proposition also relates closely to ideas in Fulton ~\cite{ful2}, Proposition 1: Use notation from the preceding discussion. Fulton shows that assuming $\prod_{j=1}^s\omega_{I^j}\neq 0$, we have the implication
$$\prod_{j=1}^s\omega_{K^j}\neq 0\in H^*(Gr(d,r)) \Rightarrow\prod_{j=1}^s\omega_{L^j}\neq 0\in H^*(Gr(d,n)).$$
But the discussion above shows that in addition, if $\prod_{j=1}^s\omega_{I^j}\neq 0$ and $\prod_{j=1}^s\omega_{K^j}\neq 0$, then the (homological) degree of the Poincar\'{e} dual of the product $\prod_{j=1}^s\omega_{L^j}$ is at least as great as the (homological) degree of the Poincar\'{e} dual of $\prod_{j=1}^s \omega_{K^j}$.
 \section{Tangent spaces and Intersection theory}
We make precise and explore the idea that if we can get a collection of Schubert varieties (with respect to different flags) in a Grassmannian $\Gr(r,W)$ to intersect transversally at a point, then the classes of these Schubert varieties have a non-zero product in the Chow ring (this criterion also appears in a paper of K. Purbhoo ~\cite{purbhoo}).

Let $V\subseteq  W$ be an $r$-dimensional subspace of an $n$-dimensional vector space $W$. Recall that the tangent space $T(\Gr(r,W))_V$ of $\Gr(r,W)$ at $V$ is canonically isomorphic to
$\home(V,W/ V)$. The following lemma (see ~\cite{sottile}, Section 2.7 for a proof) identifies the tangent spaces of Schubert cells:
\begin{lemma}\label{sottile}Let $I=\{i_1<\dots<i_r\}$  be a subset of $[n]$ of cardinality $r$ and $W$ an $n$-dimensional vector space. Let $E_{\sssize{\bullet}}$ be a complete flag on $W$ and  let $V\in\Omega^{o}_{I}(E_{\sssize{\bullet}})$. Let $E_{\bull}(V)$ and $E_{\bull}(W/ V)$ denote the induced flags on $V$ and $W/ V$ respectively. Then the tangent space to the Schubert cell $\Omega^o_{I}(E_{\bull})$ at $V$,
$$T(\Omega^o_I(E_{\bull}))_V \subseteq  T(\Gr(r,W))_V= \home(V,W/ V)$$ is given by
$$\{\phi \in \home(V,W/V)\mid \phi(E_a(V)) \subseteq  E_{i_a -a}(W/V) \text{ for }
a=1,\dots, r\}.$$
\end{lemma}
If $\me\in \Fl(W)^s$ and $I^1,\dots,I^s$ are the unique subsets of $[n]$ each of cardinality $r$ such that $V\in \cap_{j=1}^s \Omega^o_{I^j}(E^j_{\bull})$, it follows that the tangent space at $V$ to the scheme theoretic intersection $\cap_{j=1}^s\Omega^o_{I^j}(E^j_{\bull})$ is given by
\begin{equation}\label{tangent}
\{\phi \in \home(V,W/V)\mid \phi(E^j_a(V)) \subseteq  E^j_{i^j_a -a}(W/V) \text{ for }a=1,\dots,r,\text{ }j=1,\dots,s\}.
\end{equation}
\subsection{Relating tangent spaces to intersection theory}
Inspired by the above description of the tangent space of an intersection of Schubert cells (Equation ~\ref{tangent}), we make the following definition:
\begin{definition}\label{deftwo}
Let $\mi=(I^1,\dots,I^s)$ be a $s$-tuple of subsets of $[n]$ of cardinality $r$ each. Let $V$  and $Q$ be vector spaces of rank $r$ and $n-r$ respectively and $(\mf,\mg)\in \Fl(V)^s\times \Fl(Q)^s$. Define
\begin{equation}\label{intersection}
\home_{\mi}(V,Q,\mf,\mg)=\bigcap_{j=1}^s \{\phi \in \home(V,Q)\mid \phi(F^j_a) \subseteq  G^j_{i^j_a -a}\text{ for }a=1,\dots,r\}.
\end{equation}
\begin{equation}\label{vs}
=\{\phi \in \home(V,Q)\mid \phi(F^j_a) \subseteq  G^j_{i^j_a -a}\text{ for }a=1,\dots,r,\text{ }j=1,\dots,s\}
\end{equation}
\end{definition}
\begin{proposition}\label{april2}Let $V=\kappa^r$ and $Q=\kappa^{n-r}$ and $\mi=(I^1, \dots, I^s)$ a $s$-tuple of subsets of $[n]$ each of cardinality $r$. Consider the following properties:
\begin{enumerate}
\item[$(\alpha)$] $\prod_{j=1}^s\omega_{I^j}\neq 0\in H^*(Gr(r,n))$.
\item[$(\beta)$] There exists $\me\in\Fl(\kappa^n)^s$ with the following property: There is a point $T\in\bigcap_{j=1}^s \Omega^o_{I^j}(E^j_{\bull})$ such that the (smooth) varieties $\Omega^o_{I^j}(E^j_{\bull})$ for $j=1,\dots,s$ meet transversally at $T$.
\item[$(\gamma)$] If $(\mf,\mg)$ is a generic point of $\Fl(V)^s\times\Fl(Q)^s$, the rank of the vector space $\hoo$
is  $\dim(Gr(r,n))-\sum_{j=1}^s \codim(\omega_{I^j})$.
\end{enumerate}
 The implications   $(\gamma)\Rightarrow(\beta)\Rightarrow (\alpha)$ hold in any characteristic. In characteristic $0$, $(\alpha)\Rightarrow(\gamma)$.
\end{proposition}
\begin{proof}
Let $W=V\oplus Q$. We first prove $(\alpha)\Rightarrow (\gamma)$ in characteristic $0$. Choose $\me\in\Fl(W)^s$ so that $\cap_{j=1}^s\Omega^o_{I^j}(E^j_{\bull})$ is a non empty transverse intersection (using Proposition ~\ref{kleiman}). Using the group action, we may assume that
$V\in\cap_{j=1}^s\Omega^o_{I^j}(E^j_{\bull})$. Transversality and Equation ~\ref{tangent} tell us that, identifying $Q$ with $W/V$,
$$\{\phi \in \home(V,Q)\mid \phi(E^j_a(V)) \subseteq  E^j_{i^j_a -a}(Q) \text{ for }a=1,\dots,r,\text{ }j=1,\dots,s\}$$
is of rank $\dim(Gr(r,n))-\sum_{j=1}^s \codim(\omega_{I^j})$. The point $(\me(V),\me(Q))$ therefore satisfies the conditions in $(\gamma)$. It remains to see that the property in $(\gamma)$ is an open condition, i.e the set of $(\mf,\mg)$ satisfying it is an open subset of $\Fl(V)^s\times \Fl(Q)^s$. This follows from Lemma ~\ref{april1} below.

We now prove $(\gamma)\Rightarrow (\beta)$: Given $(\mf,\mg)\in \Fl(V)^s\times \Fl(Q)^s$ as in $(\beta)$, we find using Lemma ~\ref{flags}, some $\mathcal{E}\in \Fl(W)^s$  for which
\begin{enumerate}
\item[(a)] $V\in \cap_{j=1}^s\Omega^o_{I^j}(E^j_{\bull})$.
\item[(b)] The (ordered) collection of flags induced by $\mathcal{E}$ on $V$ and $Q$ are $\mf$ and $\mg$ respectively.
\end{enumerate}
The tangent space at $V$ of $\cap_{j=1}^s\Omega^o_{I^j}(E^j_{\bull})$ (see Equation ~\ref{tangent}) is then
$$\{\phi \in \home(V,Q)\mid \phi(F^j_a) \subseteq  G^j_{i^j_a -a} \text{ for }a=1,\dots,r,\text{ }j=1,\dots,s\}$$
Hence the assumption in $(\beta)$ implies that the varieties $\Omega^o_{I^1}(E^1_{\bull})$,$\dots$, $\Omega^o_{I^s}(E^s_{\bull})$ meet transversally at $V$. This proves $(\beta)$.

For $(\beta)\Rightarrow (\alpha)$. Let $\me$ and $T$ be as in $(\beta)$. Consider the intersection
$$\Omega=\bigcap_{j=1}^s \Omega_{I^j}(E^j_{\bull})$$
The assumption of transversality in $(\beta)$ tells us there is a
unique irreducible component $Z$ of $\Omega$ containing  $T$ and that $Z$ is a proper  irreducible component of the intersection $\Omega$ (see Section ~\ref{faizfaiz}). Let
$Z_1,\dots,Z_k$ be the other distinguished varieties (in the sense of
~\cite{int}, Section 6.1) of the intersection $\Omega$.

 According to Fulton-MacPherson's intersection theory (see ~\cite{int}, Section 6.1), the intersection product $\omega_{I^1}\cdot\ldots\cdot \omega_{I^s}$ can be written as a sum $m_Z[Z]+\sum_{l=1}^k \alpha_l$ where for $l=1,\dots,k$, $\alpha_l$ is a push forward of a cycle on $Z_l$. By the transversality assumption in $(\beta)$ and  ~\cite{int}, Proposition $8.2$, $m_Z=1$. By ~\cite{int}, Section $12.2$, each $\alpha_l$ is a non-negative cycle (this part uses the fact that Grassmannians are homogenous spaces and hence their tangent bundles are generated by global  sections). It is now easy to see that $\omega_{I^1}\cdot\ldots\cdot \omega_{I^s}\neq 0$ (in cohomology or the Chow ring).
\end{proof}
\begin{lemma}\label{april1} With notation from Proposition ~\ref{april2},
\begin{enumerate}
\item For $(\mf,\mg)\in \Fl(V)^s\times \Fl(Q)^s$, the following inequality holds: $$\rk(\hoo)\geq \dim(\Gr(r,n))-\sum_{j=1}^s \codim(\omega_{I^j})$$
\item The subset of   $\Fl(V)^s \times \Fl(Q)^s$ formed by $(\mathcal{F},\mathcal{G})$ such that equality holds in (1), is open (possibly empty).
\end{enumerate}
\end{lemma}
\begin{proof}Let $X=\Fl(V)^s\times\Fl(Q)^s$. Let $\mathcal{B}=\home(V,Q)\tensor\mathcal{O}_X$ and for $j=1,\dots,s$, let $\mpp^j$ be the subbundle of $\mb$ whose fiber
over a point $(\mathcal{F},\mathcal{G})$ of $X$  is
$$\{\phi\in \home(V,Q)\mid \phi(F^j_{a})\subseteq  G^j_{i^j_a-a},\text{ }a=1,\dots, r\}.$$

One verifies easily by linear algebra that $\rk(\mpp^j)=\sum_{a=1}^r (i^j_a-a).$ Now define $\mathcal{C}= \bigoplus_{j=1}^{j=s} \mathcal{B}/\mpp^j.$ Let $t:\mathcal{B}\to \mathcal{C}$  be the morphism obtained by (direct) summing the quotient morphisms $\mathcal{B}\to \mathcal{B}/\mpp^j$.

The vector space $\hoo$ is the kernel of $t$ at the point $(\mf,\mg)$, and hence is of rank at least as great as $\rk(\mathcal{B})-\sum_{j=1}^s [\rk(\mathcal{B}/\mpp^j)]$. To prove (1) it suffices to remark that $\rk(\mathcal{B})=\dim(\Gr(r,n))$ and $\rk(\mathcal{B}/\mpp^j)=\codim(\omega_{I^j})$ for $j=1,\dots, s$.

Let $U$ be the subset of $\Fl(V)^s\times\Fl(Q)^s$ appearing in the statement (2) of the lemma. We observe that $(\mf,\mg)\in U$ if and only if the map $t$ is surjective at $(\mf,\mg)\in X$. $U$ is therefore an open subset of $X$.
\end{proof}
\begin{lemma}\label{flags} Given $W$ an $n$-dimensional vector space, $V\subseteq  W$ an $r$-dimensional subspace and
$I$ a subset of $[n]$ of
cardinality $r$, let $A$ be the set  of complete flags $E_{\sssize{\bullet}}$ on $W$ such that $V \in \Omega^o_{I}(E_{\bull})$.
Then,  the natural map (a complete flag on $W$ induces a complete flag on $V$ and a complete flag on  $W/V$):
$$ \tau:A\to \Fl(V) \times \Fl(W/V)$$
is surjective.
\end{lemma}
\begin{proof}
Consider the group $P=\{g\in \GL(n)\mid gV=V\}$ which surjects to $L=\GL(V)\times \GL(W/V)$. The map $ A\to \Fl(V) \times \Fl(W/V)$ is $P$-equivariant. The action of $P$ on $\Fl(V) \times \Fl(W/V)$ is via its quotient $L$. $L$ acts transitively on $\Fl(V) \times \Fl(W/V)$ and hence the $P$-equivariance of $\tau$ shows its surjectivity.
\end{proof}
\subsection{Dimension computations}\label{dimensioncount}
The vector space $\home_{\mi}(V,Q,\mf,\mg)$ (from Definition
~\ref{deftwo}) can be  viewed as an intersection (see Equation
~\ref{intersection}) of subspaces of $\home(V,Q)$. Now,
$GL(V)\times GL(Q)$ acts with finitely many orbits (labelled by
the rank) on $\home(V,Q)$. We can obtain information on the
rank of $\hoo$, by using Proposition ~\ref{kleiman} on each
orbit. We start with a preparatory lemma:
\begin{lemma}\label{k1}
Let $V$ be an $r$-dimensional vector space and $\mf$ a generic point in $\Fl(V)^s$. Suppose $S\subseteq  V $ is a  subspace of rank $d$, and $\mk=(K^1,\dots,K^s)$ the unique $s$-tuple of subsets of $[r]$ each of cardinality $d$ such that $S\in \cap_{j=1}^s\Omega^o_{K^j}(F^j_{\bull})$. Then,
\begin{enumerate}
\item $\dim(S,V,\mf)\geq 0$.
\item $\prod_{j=1}^s  \omega_{K^j}\neq 0\in H^*(\Gr(d,V)).$
\end{enumerate}
\end{lemma}
\begin{proof}
 Genericity of $\mf$ implies that $\cap_{j=1}^s \Omega^o_{K^j}(F^j_{\bull})$ is a  proper intersection. It is non empty because $S$ is in it. Properness of intersection says that the  dimension of any irreducible component (which is non negative!) is  equal to $\dim(\Gr(d,r))-\sum_{j=1}^s \codim(\omega_{K^j})$ which is $\dim(S,V,\mf)$ by definition of the latter. This gives (1).

The genericity of $\mf$ and the nonemptiness  of $\cap_{j=1}^s \Omega^o_{K^j}(F^j_{\bull})$ give (2).
\end{proof}
The following is our main dimension count, and the proof of it will be given in Section ~\ref{countproof}. Let the notation be as in Definition ~\ref{deftwo}.
\begin{proposition}\label{propos}
 Let $(\mf,\mg)$ be a generic point in $\Fl(V)^s\times \Fl(Q)^s$. Let $\phi$ be a generic element of
$\home_{\mi}(V,Q,\mf,\mg)$. Let $S=\ker(\phi)$, $d=\rk(S)$, and let $\mk=(K^1,\dots, K^s)$ be the unique $s$-tuple of subsets of $[r]$ each of cardinality $d$ such  that  $S\in \cap_{j=1}^s\Omega^o_{K^j}(F^j_{\bull})$. Then, the rank of $\home_{\mi}(V,Q,\mf,\mg)$ is
\begin{equation}\label{dimensioncalcul}
[\dim(Gr(r,n))-\sum_{j=1}^s \codim(\omega_{I^j})]+\dim(S,V,\mf)+ \{\sum_{j=1}^s\sum_{a\in K^j}(n-r+a-i^j_a)-d(n-r)\}.
\end{equation}
\end{proposition}
The first term in Equation ~\ref{dimensioncalcul} is the ``expected rank'' of $\home_{\mi}(V,Q,\mf,\mg)$ (Lemma ~\ref{april1}). By Lemma ~\ref{k1}, $\prod_{j=1}^s \omega_{K^j}\neq 0\in H^*(\Gr(d,r))$ and $\dim(S,V,\mf)\geq 0$. The last term of Equation ~\ref{dimensioncalcul} is the same as the quantity in Inequality $(\dagger_{\mk}^{\mi})$.
\subsection{The Strategy}\label{strategem}
We explain our strategy in proving (B)$\Rightarrow$(A) in Theorem ~\ref{MainTe}. Let (B) hold in Theorem ~\ref{MainTe}. Let $V$ and $Q$ be vector spaces of ranks $r$ and $n-r$ respectively. Let $(\mf,\mg)$ be a generic point of $\Fl(V)^s\times \Fl(Q)^s$, and  $S$, $\phi$, $d$ and $\mk=(K^1,\dots,K^s)$ be as in Proposition ~\ref{propos}.

By Lemma ~\ref{k1} applied to the inclusion $S\subseteq  V$ and $\mf\in \Fl(V)^s$, we find that  $\prod_{j=1}^s  \omega_{K^j}\neq 0$. Inequality $(\dagger_{\mk}^{\mi})$: $\sum_{j=1}^s\sum_{a\in K^j}(n-r+a-i^j_a)-d(n-r)\leq 0$ therefore holds. Hence an inspection of the formula ~\ref{dimensioncalcul}, tells us that the vector space $\hoo$ is of rank not greater than 
$ [\dim(Gr(r,n))-\sum_{j=1}^s \codim(\omega_{I^j})]+\dim(S,V,\mf).$

If $\dim(S,V,\mf)=0$, then (A) holds by Lemma ~\ref{april1}, (1) and Proposition ~\ref{april2}, $(\gamma)\Rightarrow(\alpha)$.  Suppose that $\dim(S,V,\mf)>0$. Each irreducible component of the intersection $\cap_{j=1}^s\Omega^o_{K^j}(F^j)$ passing through $S$ is therefore positive dimensional. If $\psi$ is a generic element in the tangent space of $\cap_{j=1}^s\Omega^o_{K^j}(F^j)\subseteq  \Gr(d,V)$ at $S$, we can view $\psi$ as a map $S\to V/S$ (the tangent space to $\Gr(d,V)$ at $S$ is $\home(S,V/S)$) and hence  form the composite
$S\leto{\psi} V/S\leto{\phi} Q$. It turns out that the image of $\phi\circ\psi\subseteq \im(\phi)$ ``conflicts'' Proposition ~\ref{l2}, (2) with $W=Q$, $V=\im(\phi)$, $S=\im(\phi\circ \psi)$ and $\me=\mg$. This ``conflict'' is delicate,  requires a further analysis of the kernel of $\psi$ and needs to be formulated as an inductive statement. We do this in Section ~\ref{saturday} (Theorem ~\ref{stonger}) but mention an important special case here: If $\ker(\psi)=0$, let $M=\im(\phi)$ and $M^{(1)}=\im(\phi\circ \psi)\subseteq  M$. It follows from Claim ~\ref{clm} that $\dim(M^{({1})},M,\mg(M))-\dim(M^{({1})},Q,\mg)$ is of rank atleast as great as 
$$\dim(S,V,\mf)+ \{\sum_{j=1}^s \sum_{a\in K^{j}}(n-r+a-i^j_a)-d(n-r)\}$$
(where the term in the curly brackets is the same as the quantity appearing in  Inequality $(\dagger_{\mk}^{\mi})$). But since $\mg\in\Fl(Q)^s$ is  generic, we have  by Proposition ~\ref{l2}, (2), $\dim(M^{({1})},M,\mg(M))-\dim(M^{({1})},Q,\mg)\leq 0.$ Therefore,
$$\dim(S,V,\mf)+ \{\sum_{j=1}^s \sum_{a\in K^{j}}(n-r+a-i^j_a)-d(n-r)\}\leq 0.$$ The rank  of $\hoo$ is hence (see Equation ~\ref{dimensioncalcul}) not greater than
$\dim(Gr(r,n))-\sum_{j=1}^s \codim(\omega_{I^j}).$
 But by Lemma~\ref{april1}, the rank of $\hoo$ is at least as great as $\dim(Gr(r,n))-\sum_{j=1}^s \codim(\omega_{I^j}).$  Therefore, $\hoo$ is of rank $\dim(Gr(r,n))-\sum_{j=1}^s \codim(\omega_{I^j})$, and we conclude that (A) holds in Theorem ~\ref{MainTe} by 
using Proposition ~\ref{april2} $(\gamma)\Rightarrow (\alpha)$.
\subsection{Proof of Proposition ~\ref{propos}}\label{countproof}
The reader is advised to postpone the proof of proposition ~\ref{propos} and skip to Section 3 in the first reading. We prove the following statement first 
\begin{proposition}\label{strong} Consider a $6$-tuple of the form $(V,Q,\mf,\mg,\mi,\mk)$ where  $V$, $Q$, $\mf$, $\mg$ and $\mi$ are as in Definition ~\ref{deftwo} and  $\mk=(K^1,\dots,K^s)$ a $s$-tuple of subsets of $[r]$ each of cardinality $d$.
\begin{enumerate}
 \item If  $\cap_{j=1}^s\Omega^o_{K^j}(F^j_{\bull})\neq \emptyset$, then
$\rk(\hoo)$ is at least as great as $$[\dim(\Gr(r,n))-\sum_{j=1}^s \codim(\omega_{I^j})]+\{\sum_{j=1}^s\sum_{a\in K^j}(n-r+a-i^j_a)-d(n-r)\}.$$
(the element in the curly brackets is the quantity appearing in Inequality $(\dagger_{\mk}^{\mi})$.)
\item There exists a nonempty open subset $U$ of $\Fl(V)^s\times\Fl(Q)^s$ such that for $(\mf,\mg)\in U$, the dimension of each irreducible component (if non empty) of
\begin{equation}\label{latenight}
\{\phi\in\hoo\mid \rk(\ker(\phi))=d,\text{ } \ker(\phi)\in \bigcap_{j=1}^s \Omega^o_{K^j}(F^j_{\bull})\}
\end{equation}
is given by

\begin{equation}\label{sunday}
[\dim(Gr(r,n))-\sum_{j=1}^s \codim(\omega_{I^j})]
+ [\dim(\Gr(d,r))-\sum_{j=1}^s \codim(\omega_{K^j})]
\end{equation}
$$+ [\sum_{j=1}^s\sum_{a\in K^j}(n-r+a-i^j_a)-d(n-r)].$$
\end{enumerate}
\end{proposition}
{\bf Proposition ~\ref{propos} follows from Proposition ~\ref{strong}}: Let $V$, $Q$, $(\mf,\mg)$, $\mk$, $\phi$ and $S$  be as in Proposition ~\ref{propos}. Now apply Proposition ~\ref{strong}, (2) to the tuple $(V,Q,\mf,\mg,\mi,\mk)$. The subset of $\hoo$
given by Expression ~\ref{latenight} is open and dense. Hence the rank of $\hoo$ is given by Expression ~\ref{sunday} which coincides with Expression ~\ref{dimensioncalcul} (the middle term of Expression ~\ref{sunday} is in this case equal to $\dim(S,V,\mf)$).

Proposition ~\ref{strong}, together with Proposition ~\ref{april2} can be used to obtain  a different proof of (A)$\Rightarrow$ (B) in Theorem ~\ref{MainTe}. Without loss of generality we may assume that the base field is algebraically closed of  characteristic $0$. Let $V=\kappa^r$ and $Q=\kappa^{n-r}$. Assume that  $\prod_{j=1}^s\omega_{I^j}\neq 0\in H^*(\Gr(r,n))$. Let $\mk=(K^1,\dots,K^s)$, $d$ and $r$ as in Theorem ~\ref{MainTe} (B) with $\prod_{j=1}^s \omega_{K^j}\neq 0\in H^*(\Gr(d,r))$. Choose a generic point $(\mf,\mg)\in \Fl(V)^s\times
\Fl(Q)^s$. Now since $\mf$ is generic, the intersection $\cap_{j=1}^s\Omega^o_{K^j}(F^j_{\bull})\neq \emptyset$. Now apply Proposition ~\ref{strong} (1) to get that the rank of $\hoo$ is at least as great as 
\begin{equation}\label{cher}
[\dim(\Gr(r,n))-\sum_{j=1}^s \codim(\omega_{I^j})]+\sum_{j=1}^s\sum_{a\in K^j}(n-r+a-i^j_a)-d(n-r).
\end{equation}
Since $\prod_{j=1}^s\omega_{I^j}\neq 0\in H^*(\Gr(r,n))$, condition $(\alpha)$ of Proposition ~\ref{april2} holds and hence condition $(\gamma)$ holds as well. Therefore, the  rank of $\hoo$ is
$[\dim(Gr(r,n))-\sum_{j=1}^s \codim(\omega_{I^j})]$
and a comparison with Inequality ~\ref{cher} yields Inequality ($\dagger_{\mk}^{\mi}$).

We now proceed to the proof of Proposition ~\ref{strong} 
\begin{proof}
Let $\ms$ be the universal subbundle of $\mv=V\tensor \mathcal{O}_{X}$
on  $X=\Gr(d,V)$. Let $Y$ be the total space of the vector bundle $\shom(\mv/\ms,Q\tensor\mathcal{O}_X)$ over $X$. Let $Y^o$ be the open subset of $Y$ formed by points where the universal morphism $(\mv/\ms)_Y\to Q_Y$ is an injection (on the fibers). We can think of points in $Y$ as pairs $(T,\psi)$ where
$T\subseteq  V$ is a $d$-dimensional subspace and $\psi:V\to Q$ is a homomorphism satisfying $\psi(T)=0$. The points in $Y^o$ are pairs $(T,\psi)\in Y$ with $\ker(\psi)=T$. It is easy to see that $\dim(Y)=\dim(X)+(n-r)(r-d)$. $\GL(V)\times \GL(Q)$ acts on $Y$ and this action is transitive when restricted to $Y^o$.
For $j=1,\dots,s$, let
$P^j(F^j_{\bull},G^j_{\bull})$ be the subset of $Y$ formed by elements $(T,\psi)$ satisfying
\begin{enumerate}
\item[(i)]$T\in \Omega^o_{K^{j}}(F^{j}_{\bull}).$
\item [(ii)]$\psi(F^j_a)\subseteq  G^j_{i^j_a-a}$ for $a=1,\dots,r.$
\end{enumerate}
It is easy to see that $P^j(F^j_{\bull},G^j_{\bull})\neq \emptyset$ (because $\psi=0$ is an acceptable choice of $\psi$  in (ii) above). Let $\delta^j:P^j(F^j_{\bull},G^j_{\bull})\to \Omega^o_{K^j}(F^j_{\bull})$. This map is equivariant with respect to the Borel subgroup of $\GL(V)$ determined by $F^j_{\bull}$ and is hence surjective. We compute the dimensions of the fibers of this morphism. Let $j$ be fixed for this calculation. If $T\in\Omega^o_{K^j}(F^j_{\bull})$, pick elements $e(a)\in F^j_a$ for $a=1,\dots,r$ such that
\begin{enumerate}
\item [(a)]$F^j_a=\Spann\{e(1),\dots,e(a)\}$ for $a=1,\dots,r$.
\item [(b)]For $t=1,\dots,d$, $e(k^j_t)\in T\cap F^j_{k^j_t}$.
\end{enumerate}
 If $\psi\in(\delta^j)^{-1}(T)$, then $\psi$  has to send $e(k^j_t)$ to $0$ and for $a\in [r]\smallsetminus K^j$, $\psi(e(a))$ should belong to $G^j_{{i^j_a-a}}$. So the fiber of $\delta^j$ over $T$ is the vector space $\prod_{a\in [r]\smallsetminus K^j} G^j_{{i^j_a-a}}$. It is now easy to see that $P^j(F^j_{\bull},G^j_{\bull})$ is a smooth and irreducible subvariety  of $Y$ of dimension
\begin{equation}\label{dawn}
\dim(\Omega_{K^j}(F_{\bull}))+ \sum_{a\in[r]\smallsetminus  K} (i^j_a-a).
\end{equation}
With $S$ as in (1), consider the vector subpace of $\hoo$
$$A=\{\phi\in \hoo\mid \phi(S)=0\}\subseteq  \home(V/S,Q)$$
Clearly $A=\cap_{j=1}^s(\delta^j)^{-1}(S)$. The rank of the vector space $A$ is therefore at least as great as
$$\rk(\home(V/S,Q))-\sum_{j=1}^s(\rk(\home(V/S,Q))-\rk((\delta^j)^{-1}(S)))$$
which we rewrite as
\begin{equation}\label{l011}
(n-r)(r-d)-\sum_{j=1}^s[(n-r)(r-d)-\sum_{a\in[r]\smallsetminus  K^{j}}{(i^j_a-a)}]
\end{equation}
$$=(n-r)(r-d)-\sum_{j=1}^s\sum_{a\in[r]\smallsetminus  K^{j}}(n-r+a-i^j_a)$$
$$=\{r(n-r)-\sum_{j=1}^s\sum_{a=1}^r(n-r+a-i^j_a)\}+\{-d(n-r)+\sum_{j=1}^s\sum_{a\in K^{j}}(n-r+a-i^j_a)\}.$$
Clearly $\rk(\hoo)\geq\rk(A)$ and this proves (1).

For (2), denote by $H$ the subspace of $\hoo$ given by Expression ~\ref{latenight}. Clearly, $\hmk=\cap_{j=1}^s P^j(F^j_{\bull},G^j_{\bull})\cap Y^o$. To compute the dimension of irreducible components of $\hmk$, we observe first that for $j=1,\dots,s$, (see Equation ~\ref{dawn})
$$\codim(P^j(F^j_{\bull},G^j_{\bull}),Y)=\dim(Y)-\dim(\Omega_{K^j}(F^j_{\bull}))-\sum_{a\in[r]\smallsetminus  K^j}{(i^j_a-a)}.$$
$$=d(r-d) +(n-r)(r-d)-\dim(\Omega_{K^j}(F^j_{\bull}))-\sum_{a\in[r]\smallsetminus  K^j}{(i^j_a-a)}$$
$$=\codim(\omega_{K^j})+\sum_{a\in [r]\smallsetminus K^j}(n-r-{(i^j_a-a)}).$$

Now, $\GL(V)\times\GL(Q)$ acts transitively on $Y^o$, therefore the genericity of $(\mf,\mg)$ and Proposition ~\ref{kleiman} imply
 that each irreducible component of $\hmk$ is of dimension
$$[\dim(X)+(n-r)(r-d)]-\sum_{j=1}^s\codim(\omega_{K^j})- \sum_{j=1}^s\sum_{a\in[r]\smallsetminus  K^j}((n-r)-{(i^j_a-a)})$$
$$=[\dim(\Gr(d,r))-\sum_{j=1}^s\codim(\omega_{K^j})]+ \{(n-r)(r-d)-\sum_{j=1}^s\sum_{a\in[r]\smallsetminus  K^j}((n-r)-{(i^j_a-a)})\}.$$
We now conclude the proof using Equation ~\ref{l011}.
\end{proof}
\section{Formulation of the main result}\label{saturday}
In this section we state our main theorem and see how it implies (B)$\Rightarrow$(A) in Theorem ~\ref{MainTe}.
\begin{theorem}\label{stonger} Consider a $5$-tuple of the form $(V,Q,\mf,\mg,\mi)$ where $V$ and $Q$ are non-zero vector spaces of ranks $r$ and $n-r$ respectively, $\mi=(I^1,\dots,I^s)$ a $s$-tuple of subsets of $[n]$ each of cardinality $r$ and   $(\mf,\mg)\in \Fl(V)^s\times \Fl(Q)^s$ a generic point (see Section ~\ref{strategem2}).

Then, there exists a   filtration by vector subspaces
\begin{equation}\label{filtratione}
S^{(h)}\subsetneq S^{(h-1)}\subsetneq\dots \subsetneq S^{(1)}\subsetneq S^{(0)}=V
\end{equation}
and injections (of vector spaces) from the graded quotients $\eta_{u}:{S^{(u)}}/{S^{(u+1)}}\hookrightarrow Q$ for $u=0,\dots, h-1$ such that the following property is satisfied: For $u=1,\dots,h$, let $d_u$ be the rank of $S^{(u)}$, ${\Ftwo}(u)=({\Fone}^1(u),\dots,{\Fone}^s(u))$ the unique $s$-tuple of subsets of $[r]$ each of cardinality $d_u$ such that  $S^{(u)}\in \cap_{j=1}^s\Omega^o_{{\Fone}^j(u)}(F^j_{\bull})\subseteq\Gr(d_u,V)$, then
\begin{enumerate}
\item[(i)] $\dim(S^{(h)},V,\mf)=0.$
\item[(ii)] For $u=0,\dots,h-1$, $j=1,\dots,s$ and $a=1,\dots,r$, (where we write $\eta_u$ again  for the composite $S^{(u)}\to {S^{(u)}}/{S^{(u+1)}}\leto{\eta_u} Q$)
$$\eta_u(S^{(u)}\cap F^j_a)\subseteq  G^j_{i^j_a-a}$$
\item[(iii)] The vector space $\hoo$ is of rank (the second term of the expression below is the same as the quantity appearing in Inequality $(\dagger_{{\Ftwo}(h)}^{\mi})$)
\begin{equation}\label{une}
[\dim(Gr(r,n))-\sum_{j=1}^s \codim(\omega_{I^j})]+ \{\sum_{j=1}^s\sum_{a\in {\Fone}^j(h)}(n-r+a-i^j_a)-d_h(n-r)\}.
\end{equation}
\end{enumerate}
\end{theorem}
\begin{remark} In Theorem ~\ref{stonger}, the case $h=0$ can occur. Also, $S^{(h)}$ can be the $0$ vector space (see Example ~\ref{fourr}). The filtration is not necessarily unique.
\end{remark}
\subsection{Theorem ~\ref{stonger} implies  (B)$\Rightarrow$ (A) in Theorem ~\ref{MainTe}}\label{second}
Let $\mi=(I^1,\dots,I^s)$ be a $s$-tuple of subsets of $[n]$ each of cardinality $r$ as in   Theorem ~\ref{MainTe}. Further assume that Condition (B) holds in Theorem ~\ref{MainTe}.

Let $V$ and $Q$ be vector spaces of ranks $r$ and $n-r$ respectively and $(\mf,\mg)$ a generic point of $\Fl(V)^s\times\Fl(Q)^s$. Apply Theorem ~\ref{stonger} to the $5$-tuple $(V,Q,\mf,\mg,\mi)$ and use the same notation as in the conclusion of Theorem ~\ref{stonger}. The rank of the vector space $\hoo$ is (by conclusion (iii) of Theorem ~\ref{stonger})
$$[\dim(Gr(r,n))-\sum_{j=1}^s \codim(\omega_{I^j})]+ \{\sum_{j=1}^s\sum_{a\in {\Fone}^j(h)}(n-r+a-i^j_a)-d_h(n-r)\}.$$
Lemma ~\ref{k1} (applied to $S^{(h)}\subseteq  V$ and $\mf\in\Fl(V)^s$) implies that $\prod_{j=1}^s \omega_{{\Fone}^j(h)}\neq 0\in H^*(Gr(d_h,V))$. 
By Inequality $(\dagger_{{\Ftwo}(h)}^{\mi})$), $\sum_{j=1}^s\sum_{a\in {\Fone}^j(h)}(n-r+a-i^j_a)-d_h(n-r)\leq 0$. The rank of $\hoo$ is therefore not greater than 
$[\dim(Gr(r,n))-\sum_{j=1}^s \codim(\omega_{I^j})].$ The reverse inequality  holds by Lemma ~\ref{april1}, (1) and hence the  rank of $\hoo$ equals $[\dim(Gr(r,n))-\sum_{j=1}^s \codim(\omega_{I^j})].$ We  use Proposition ~\ref{april2}, $(\gamma)\Rightarrow(\alpha)$ to conclude that  (A) holds in Theorem ~\ref{MainTe}. 
\section{First discussion of Genericity}\label{strategem2}
The proof of Theorem ~\ref{stonger} which will be given in Section ~\ref{Five}
is by induction on $r=\rk(V)$. Let $(\mf,\mg)$ be a generic point of $\Fl(V)^s\times \Fl(Q)^s$ and $\phi$ a generic element of $\hoo$. Let $S$ and $\mk$ be as in the conclusion of Proposition ~\ref{propos}. In the proof of Theorem ~\ref{stonger}, we are going to set $S^{(1)}=S$ and apply the induction hypothesis to the $5$-tuple $(S,V/S,\mf(S),\mf(V/S),\mk)$ if $0<\rk(S)<r$. But the induction hypothesis was only for generic pairs of flags in $\Fl(S)^s\times \Fl(V/S)^s$. Therefore, to justify the  above application of induction hypothesis, we need to make the notion of genericity required in Theorem ~\ref{stonger} precise, and also prove that in Proposition ~\ref{propos}, the induced pair $(\mf(S),\mf(V/S))\in \Fl(S)^s\times\Fl(V/S)^s$ is ``generic''. We achieve this in the following way:
\begin{enumerate}
\item In Section ~\ref{whatdreamsmaycome} (which the reader may want to skim over now), we show that there is a simultaneous choice of non empty open subsets $A(R,T)$ of $\Fl(R)^s\times\Fl(T)^s$ for every pair of nonzero vector spaces $(R,T)$ such that the following property is satisfied: Suppose $V$, $Q$, $\mi$, $\mf$ and  $\mg$ are as in Definition ~\ref{deftwo} with the additional condition that $(\mf,\mg)\in A(V,Q)$. Then (as in  Proposition ~\ref{propos}), if $\phi$ is a generic element of $\hoo$ and  $S=\ker(\phi)$, the induced pair $(\mf(S),\mf(V/S))$ is a point in $A(S,V/S)$. We will also require some other genericity properties which will be listed below. The main idea behind the genericity of induced structures appears in Section ~\ref{hamlett}. The construction of $A(R,T)$ is independent of the proof of Theorem ~\ref{stonger}.
\item We prove Theorem ~\ref{stonger} for any $(\mf,\mg)\in A(V,Q)\subseteq\Fl(V)^s\times\Fl(Q)^s$ (in Section ~\ref{Five}) using only the properties of $A(R,T)$ listed in Section ~\ref{ray}.
\end{enumerate}
On a first reading we suggest that the reader assume that such a (simultaneous) choice of $A(R,T)$ exists. 

In Section ~\ref{ray} we  state the properties that the open subsets $A(R,T)\subseteq \Fl(R)^s\times \Fl(T)^s$ satisfy. To do this we first list some  genericity properties and give names to the corresponding open subsets of parameter spaces.

For a vector space  $W$ of rank $n$,  define  $B(W)\subseteq  \Fl(W)^s$ to be the largest Zariski open subset of $\Fl(W)^s$ satisfying the following property: If $\me\in B(W)$ and $\mi=(I^1,\dots,I^s)$ a $s$-tuple of subsets of $[n]$ each of the same cardinality $r$, then every irreducible component of the  intersection $\cap_{j=1}^s\Omega_{I^j}(E^j_{\bull})$ (which is possibly empty) is proper. By Proposition ~\ref{kleiman}, it follows that $B(W)$ is nonempty. An element $\me \in B(W)$ is said to be {\bf generic for intersection theory} in $W$ . It is easy to see that in  Proposition ~\ref{l2}, (2) (resp. in Lemma ~\ref{k1}) the conclusions hold if $\me$ (resp. $\mf$) is generic for intersection theory in $W$ (resp. $V$).

 The proof of Proposition ~\ref{propos} produces various open subsets in the parameter spaces:  Let the notation be as in Definition ~\ref{deftwo}. Let $r(\mi)$ be the rank of  $\home_{\mi}(V,Q,\mf,\mg)$ for generic $(\mf,\mg)\in \Fl(V)^s\times\Fl(Q)^s$. Proposition ~\ref{propos} gives a nonempty open subset $O(V,Q,\mi)\subseteq \Fl(V)^s\times\Fl(Q)^s$, an integer $d\in\{0,1,\dots,r\}$, a $s$-tuple $\mk(\mi)=\mk=(K^1,\dots,K^s)$ of subsets of $[r]$ each of cardinality $d$ such that for $(\mf,\mg)\in O(V,Q,\mi)$, if $\phi$ is a generic element of $\hoo$ and $S=\ker(\phi)$, then $\rk(S)=d$ and
\begin{enumerate}
\item  $S\in \cap_{j=1}^s \Omega^o_{K^j}(F^j_{\bull}(V))$.
\item  $\rk(\home_{\mi}(V,Q,\mf,\mg))=r(\mi)$.
\item $r(\mi)$ is given by (same as  Expression ~\ref{dimensioncalcul})
$$[\dim(Gr(r,n))-\sum_{j=1}^s \codim(\omega_{I^j})]+\dim(S,V,\mf)+ [\sum_{j=1}^s\sum_{a\in K^j}(n-r+a-i^j_a)-d(n-r)].$$
\end{enumerate}
\subsection{Properties of $A(R,T)$}\label{ray}
We now list the properties of $A(R,T)$ (which will follow from their  construction in Section ~\ref{whatdreamsmaycome}):
\begin{itemize}
 \item The choice is ``functorial for isomorphisms''. That is, if $R\to R'$ and $T\to T'$ are isomorphisms of vector spaces, then the induced bijection $\Fl(R)^s\times\Fl(T)^s\to \Fl(R')^s\times \Fl(T')^s$ carries $A(R,T)$ bijectively to $A(R',T')$. This functorial property clearly implies that $A(R,T)$ is a $\GL(R)\times\GL(T)$ stable subset of $\Fl(R)^s\times\Fl(T)^s$.
\item Let $V$ and $Q$ be  vector spaces of arbitrary ranks $r$ and $m$ respectively and $(\mf,\mg)\in A(V,Q)$. Then,
\begin{enumerate}
\item[(G1)] $\mf$ is generic for intersection theory in $V$, $\mg$ is generic for intersection theory in $Q$.
\item[(G2)] For any choice of $s$-tuple of subsets $\mi=(I^1,\dots,I^s)$ of $[r+m]$ each of cardinality $r$, $(\mf,\mg)\in O(V,Q,\mi)$.
\item[(G3)] For $\mi$ as in (G2), there is a nonempty  open subset of $\hoo$ such that for $\phi$ in this open subset, if $0<\rk(\ker(\phi))<r$, then letting   $S=\ker(\phi)$, we have $(\mf(S),\mf(V/S))\in A(S,V/S)$.
\end{enumerate}
\end{itemize}
\section{Proof of Theorem ~\ref{stonger}}\label{Five}
Proof by induction on $r$. We will prove the theorem for any $(\mf,\mg)\in A(V,Q)$ (whose properties are described in Section ~\ref{strategem2}). Assume that we have proved this result for all values of $r$ satisfying $0<r<r_0$ and prove it for $r=r_0$. The base case for induction is obtained by setting $r_0=1$ in the argument below. A few examples are given in  Section ~\ref{ILLUSTRIOUS} to illustrate the constructions in the proof.

If $\rk(\hoo)=0$, the filtration is just the singleton $V$ and $h=0$, so no maps $\eta$ need to be given. Clearly, the condition in (iii) is met (the two bracketed terms in Equation ~\ref{une} cancel each other in this case).

Now, assume that $\rk(\hoo)\neq 0$. Let $\phi\in \hoo$ be generic. Let $S=\ker(\phi)$, $d=\rk(S)$ and $\mk=(K^1,\dots,K^s)$ be the unique $s$-tuple of subsets of $[r]$ each of cardinality $d$ such that $S\in\cap_{j=1}^s \Omega^o_{K^j}(F^j_{\bull})$. Since $\mf$ is generic for intersection theory in $V$ (see Section ~\ref{strategem2}), by Lemma ~\ref{k1},
\begin{equation}\label{puccini}
\prod_{j}\omega_{K^j}\neq 0 \text{ and }  \dim(S,V,\mf)\geq 0.
\end{equation}
  If $\dim(S,V,\mf)=0$  we can take $S\subsetneq V$ to be the filtration and  $\eta_0=\phi:V/S \hookrightarrow Q$ (see Examples ~\ref{onee} and ~\ref{twoo}). This satisfies (ii) because $\phi\in\hoo$. (i) is true because $\dim(S,V,\mf)=0$. The equality required in (iii) follows from Proposition ~\ref{propos} which tells us that  the rank of $\hoo$ is ($\dim(S,V,\mf)=0$)
$[\dim(Gr(r,n))-\sum_{j=1}^s \codim(\omega_{I^j})]+ \{\sum_{j=1}^s\sum_{a\in K^j}(n-r+a-i^j_a)-d(n-r)\}.$

We now consider the case $\dim(S,V,\mf)>0$ (as in Example ~\ref{fourr}). In this case, it is easy to see that  $0<d<r$. Now by our discussion of genericity (Section ~\ref{strategem2}) the pair $(\mf(S),\mf(V/S))$ is in $A(S,V/S)$. {\bf The induction hypothesis can therefore be applied to  the $5$-tuple $(S,V/S,\mf(S),\mf(V/S),\mk)$}.

 We therefore find a filtration
$$S^{(h)}\subsetneq S^{(h-1)}\subsetneq\dots\subsetneq S^{(1)}=S$$
 and morphisms $\gamma_u:S^{(u)}/S^{(u+1)}\hookrightarrow V/S$ for $u=1,\dots, h-1$ which satisfy the conclusions of Theorem ~\ref{stonger} (we have shifted the numbering by one).

We  claim that the filtration $$S^{(h)}\subsetneq S^{(h-1)}\subsetneq\dots\subsetneq S^{(1)}=S\subsetneq S^{(0)}=V$$ and maps $\eta_u= \phi\circ\gamma_u$ for $u=1,\dots,s$, $\eta_0=\phi$ satisfy the conditions (i), (ii) and (iii) in Theorem ~\ref{stonger}. For $u=1,\dots,h$, let $d_u$ be the rank of $S^{(u)}$.
$\eta_u$ is a composite of two injections ($\gamma_u$ and $\phi:V/S\hookrightarrow Q$) and is hence an injective map.
\begin{notation}\label{notate1}{\it
\begin{enumerate}
\item Let  $\ml(u)=(L^1(u),\dots,L^s(u))$ be the unique $s$-tuple of subsets of $[d]$ each of cardinality $d_u$, such that $S^{(u)}\in \cap_{j=1}^s\Omega^o_{L^j(u)}(F^j_{\bull}(S))\subseteq  \Gr(d_u,S)$.
\item Let ${\Ftwo}(u)=({\Fone}^1(u),\dots,{\Fone}^s(u))$ be the unique $s$-tuple of subsets of $[r]$ each of cardinality $d_u$ such that  $S^{(u)}\in \cap_{j=1}^s \Omega^o_{{\Fone}^j(u)}(F^j_{\bull})\subseteq  \Gr(d_u,V)$ (${\Ftwo}(1)$ is same as $\mk$ defined above).
\end{enumerate}}
\end{notation}
 By Lemma ~\ref{falcon} (applied to $S^{(u)}\subseteq S\subseteq V$ and $\mf\in \Fl(V)^s$)
\begin{equation}\label{r1emember}
{\Fone}^j(u)=\{k^j_b\mid b\in L^j(u)\}
\end{equation}
{\bf Verification of (i)}:
We know that $\prod_{j=1}^s \omega_{K^j}\neq 0\in H^*(\Gr(d,r))$ (by Equation ~\ref{puccini}) and  the  induced collection $\mf(S)$ is generic for intersection theory in $S$ ($(\mf(S),\mf(V/S))\in A(S,V/S)$). Therefore, by Proposition ~\ref{l2} (2), Inequality ($\dagger_{\ml(h)}^{\mk}$) holds.

Inductive conclusion (iii) for the $5$-tuple $(S,V/S,\mf(S),\mf(V/S),\mk)$ tells us that the rank of $\home_{\mk}(S,V/S,\mf(S),\mf(V/S))$ is
$$\dim(\Gr(d,r))-\sum_{j=1}^s\codim(\omega_{K^j})+\{\sum_{j=1}^s \sum_{b\in L^j(h)}(r-d+b-k^j_b)-d_h(r-d)\}$$
where the term in the curly brackets is the quantity appearing in Inequality $(\dagger_{\ml(h)}^{\mk})$ and is hence $\leq 0$. On the other hand, according to Lemma ~\ref{april1} (1), the rank of $\home_{\mk}(S,V/S,\mf(S),\mf(V/S))$ is $\geq\dim(\Gr(d,r))-\sum_{j=1}^s\codim(\omega_{K^j})$. We therefore conclude that equality holds in Inequality $(\dagger_{\ml(h)}^{\mk})$:
\begin{equation}\label{E1}
\sum_{j=1}^s\sum_{b\in L^j(h)}(r-d+b-k^j_b)-d_h(r-d)=0.
\end{equation}
By Lemma ~\ref{l2}, (1), Equation ~\ref{E1} this gives  
\begin{equation}\label{E2}
\dim(S^{(h)},S,\mf(S))-\dim(S^{(h)},V,\mf)=0.
\end{equation}
 Now, by inductive conclusion (i) for the $5$-tuple $(S,V/S,\mf(S),\mf(V/S),\mk)$, we have 
$\dim(S^{(h)},S,\mf(S))=0.$
 We therefore obtain the desired conclusion: $\dim(S^{(h)},V,\mf)=0$.

\noindent{\bf Verification of (ii):} We need to verify that for ${u}=0,\dots,h-1$, $j=1,\dots,s$ and $a=1,\dots,r$,
\begin{equation}\label{dag}
\eta_{u}(S^{({u})}\cap F^j_a)\subseteq  G^j_{i^j_a-a}.
\end{equation}
Now suppose ${u}$, $j$ and $a$ are as above. If ${u}=0$, Inclusion ~\ref{dag}
is clear because $\phi\in \hoo$. So assume ${u}>0$ and find $t$ such that
$F^j_a\cap  S=F^j_t(S)$. Clearly, $k^j_t\leq a$.

From $\eta_{{u}}=\phi\circ \gamma_{u}$, we see that
$$\eta_{{u}}(S^{({u})}\cap F^j_a)=\phi\gamma_{u}(S^{({u})}\cap F^j_t(S))\subseteq  \phi (F^j_{k^j_t-t}(V/S))$$
where in the last inclusion, we have used the property (ii) satisfied by the maps $\gamma_{u}$: $\gamma_{u}(S^{({u})}\cap F^j_t(S))\subseteq  F^j_{k^j_t-t}(V/S).$ Inclusion ~\ref{dag} now follows from
$$\phi (F^j_{k^j_t-t}(V/S))=\phi((F^j_{k^j_t}+S)/S)=\phi(F^j_{k^j_t})\subseteq \phi(F^j_a)\subseteq  G^j_{i^j_a-a}.$$
{\bf Verification  of (iii):} For $u=0,\dots,h-1$, let $M^{({u})}$ be the image of $\eta_{u}$ and let $M=\text{image} (\phi) (=M^{(0)})$.
\begin{clm}\label{clm} Let
\begin{equation}\label{formule}
b({u})=\dim(S^{(u)},V,\mf)- \dim(S^{(u)},S,\mf(S)) +\{\sum_{j=1}^s\sum_{a\in {\Fone}^j(u)}(n-r+a-i^j_a)-d_u(n-r)\}
\end{equation}
$\text{for }{u}=1,\dots,h-1$ (notice that in Equation ~\ref{formule}, the term in the curly brackets is the quantity that appears in Inequality $(\dagger_{{\Ftwo}(u)}^{\mi})$).
Then, for such $u$,
$$\dim(M^{({u})},M,\mg(M))-\dim(M^{({u})},Q,\mg)\geq b({u})-b({u}+1).$$
\end{clm}
{\bf Claim ~\ref{clm} implies property (iii):}
 Now, $\mg\in \Fl(Q)^s$ is generic for intersection theory in $Q$ (see Section ~\ref{strategem2}), therefore, by  Proposition ~\ref{l2}, (2) we have $\dim(M^{({u})},M,\mg(M))-\dim(M^{({u})},Q,\mg)\leq 0$. Hence, the claim tells us that for $u=1,\dots h-1$, $b(u)\leq b(u+1)$ and therefore $b(h)\geq b(1)$.

But $\dim(S,S,\mf(S))=0$ and Proposition ~\ref{propos} tells us that $\rk(\hoo)$ equals
$b(1) +\dim(Gr(r,n))-\sum_{j=1}^s \codim(\omega_{I^j})$. On the other hand,  Equation ~\ref{E2} gives $$b(h)=\sum_{j=1}^s\sum_{a\in {\Fone}^j(h)}(n-r+a-i^j_a)-d_h(n-r).$$
Therefore, it follows from $b(h)\geq b(1)$ that  the rank of $\hoo$ is not greater than Expression ~\ref{une}. The reverse inequality holds by Proposition ~\ref{strong}, (1) (applied to $\mj(h)$). Property (iii) therefore holds.
\begin{proof}(Of the claim)
Fix a $u$ and consider
$$\frac{S^{(u)}}{S^{(u+1)}}\letof{\eta_u,\sim} M^{(u)}\subseteq  M\subseteq Q $$
The claim will be proved by applying  Lemma ~\ref{l2}, (1) to $M^{(u)}\subseteq  M  \subseteq  Q$ and $\mg\in\Fl(Q)^s$. We separate out an important part of the computation:
\begin{lemma}\label{lmma}
$\dim(M^{(u)},M,\mg(M))-\dim(M^{(u)},Q,\mg)$ is at least as great as
\begin{equation}\label{NW}
\sum_{j=1}^s\sum_{t\in L^j(u)\smallsetminus L^j(u+1)}[n-r-(r-d)+ (k^j_t-t)-(i^j_{k^j_t}-k^j_t)]-(d_u-d_{u+1})(n-r-(r-d))
\end{equation}
\end{lemma}
\begin{proof}
We introduce the following notation for the computation
\begin{itemize}
\item Let  $H^1,\dots,H^s$ be the subsets of $[n-r]$ each of cardinality $r-d$ such that $M\in \cap_{j=1}^s\Omega^o_{H^j}(G^j_{\bull})\subseteq \Gr(r-d,Q)$. 
\item  Let  ${\capn}^1,\dots,{\capn}^s$ be subsets of $[r-d]$ each of cardinality $d_u-d_{u+1}$ such that $M^{(u)}\in \cap_{j=1}^s\Omega^o_{{\capn}^j}(G^j_{\bull}(M))\subseteq \Gr(d_u-d_{u+1},M)$.
\item Let
\begin{equation}\label{yyy}
\{y^j(1)<\dots<y^j(d_u-d_{u+1})\}=L^j(u)-L^j(u+1).
\end{equation}
\end{itemize}
By Lemma ~\ref{l2} (1),
\begin{equation}\label{tvguide}
\dim(M^{(u)},M,\mg(M))-\dim(M^{(u)},Q,\mg)
\end{equation}
$$=\sum_{j=1}^s\sum_{l=1}^{d_u-d_{u+1}}[n-r-(r-d)+ p^j_{\ell}- h^j_{p^j_{\ell}}]-(d_u-d_{u+1})(n-r-(r-d)).$$

Comparing this with the required Inequality ~\ref{NW}, we need to show that for (see Equation ~\ref{yyy})
$\ell\in\{1,\dots, d_u-d_{u+1}\}$, $j\in\{1,\dots, s\}$ and $t=y^j(\ell)$, the following inequality holds
$$(n-r)-(r-d)+{\smalln}^j_{\ell}-h^j_{{\smalln}^j_{\ell}}\geq (n-r)-(r-d)+(k^j_{t}- t) -(i^j_{k^j_{t}}- {k^j_{t}})$$
Or that,
\begin{equation}\label{novemberr}
{\smalln}^j_{\ell}-h^j_{{\smalln}^j_{\ell}}\geq (k^j_{t}- t) -(i^j_{k^j_{t}}- {k^j_{t}}).
\end{equation}
For this set
\begin{itemize}
\item $a= k^j_t=k^j_{y^j(\ell)}$.
\item $x=\rk(M\cap G^j_{i^j_a-a})$ so that $M\cap G^j_{i^j_a-a}=G^j_x(M)$.
\end{itemize}
From $t=y^j(\ell)$ we find that
$\rk(F^j_{t}(S)\cap S^{(u)})-\rk(F^j_{t}(S)\cap S^{(u+1)})=\ell.$
Now, $\eta_u:\frac{S^{(u)}}{S^{(u+1)}}\leto{\sim} M^{(u)}$. Therefore $\eta_u(S^{(u)}\cap F^j_{t}(S))\subseteq M^{(u)}$ is of rank $\ell$. By property (ii) (see Inclusion ~\ref{dag}),
$$\eta_u(S^{(u)}\cap F^j_{t}(S))\subseteq \eta_u(S^{(u)}\cap F^j_a)\subseteq  M\cap G^j_{i^j_a-a}=G^j_x(M).$$
This gives  that ${\smalln}^j_{\ell}\leq x$. Now, the numbers $b-h^j_b$ weakly decrease as $b$ increases with $j$ fixed and therefore from ${\smalln}^j_{\ell}\leq x$ we obtain
\begin{equation}\label{scrooge1}
{\smalln}^j_{\ell}-h^j_{{\smalln}^j_{\ell}}\geq x-h^j_x.
\end{equation}
From the definition of $x$, $h^j_x\leq i^j_a -a$. Also, $M\cap G^{j}_{i^j_{a}-a}$ is of rank atleast as great as $k^j_t-t$ because it contains  
$\phi(F^j_{k^j_t})=\phi(F^j_{k^j_t-t}(V/S))$ which is of rank $k^j_t-t$ ($\phi(F^j_{k^j_t})\subseteq  G^{j}_{i^j_{a}-a}$ since $\phi\in \hoo$ and $a=k^j_t$). 
This gives $x\geq k^j_t-t$. We therefore obtain
\begin{equation}\label{scrooge2}
x-h^j_x \geq (k^j_{t}- t) -(i^j_{k^j_{t}}- {k^j_{t}}).
\end{equation}
 Inequalities ~\ref{scrooge1} and ~\ref{scrooge2} give  Inequality ~\ref{novemberr} and the proof of Lemma ~\ref{lmma} is complete.
\end{proof}
We write the right hand side of Expression ~\ref{NW} as $A(u)-A(u+1)$ where
$$A(u)=\sum_{j=1}^s\sum_{t\in L^j(u)}[n-r-(r-d)+ (k^j_t-t)-(i^j_{k^j_t}- k^j_t)]-d_u(n-r-(r-d)).$$
We  rearrange: $A(u)$ equals
\begin{equation}\label{height}[-d_u(n-r)+ \sum_{j=1}^s\sum_{t\in L^j(u)}(n-r+k^j_t-i^j_{k^j_t})]+\{d_u(r-d)-\sum_{j=1}^s\sum_{t\in L^j(u)}(r-d+t -k^j_t)\}.
\end{equation}
The term in the square brackets in  Equation ~\ref{height} is, using Equation ~\ref{r1emember}
$$-d_u(n-r) +\sum_{j=1}^s\sum_{a\in J^j(u)}(n-r+a-i^j_a)$$
 The term in the curly brackets in Equation ~\ref{height} is, by  Lemma ~\ref{l2}, (1)
$$\dim(S^{(u)},V,\mf)-\dim(S^{(u)},S,\mf(S)).$$
Therefore, $A(u)$ equals
$$\dim(S^{(u)},V,\mf)-\dim(S^{(u)},S,\mf(S))-d_u(n-r)+ \sum_{j=1}^s\sum_{a\in J^j(u)}(n-r+a-i^j_a)$$
which is the same as $b(u)$ (where $b(u)$ is defined in Equation ~\ref{formule}) and hence,
\begin{equation}\label{au}
A(u)-A(u+1)=b(u)-b(u+1)
\end{equation}
It now follows from Lemma ~\ref{lmma} and  Equation ~\ref{au} that
$$b(u)-b(u+1)\leq \dim(M^{(u)},M,\mg(M))-\dim(M^{(u)},Q,\mg)$$
which proves the inequality in the claim.
\end{proof}
\subsection{Illustrative examples}\label{ILLUSTRIOUS}
Unravelling the induction in Theorem ~\ref{stonger}, we find a construction of Filtration ~\ref{filtratione} with the required properties. We give some examples of this filtration in this section. The first example is one in which  $\mi$ in Theorem ~\ref{MainTe} fails to satisfy condition (A) and Theorem ~\ref{stonger} detects the failure of an  inequality $(\dagger_{\mk}^{\mi})$. 
\begin{example}\label{onee}
Consider the case $s=2$, $r=2$, $n=4$, $I^1=\{1,4\}$ and $I^2=\{2,3\}$. Let $V$ and $Q$ be vector spaces each of rank $2$ and $(\mf,\mg)$ a generic point in $\Fl(V)^2\times\Fl(Q)^2$. A map $\phi\in \home(V,Q)$ is in $\hoo$ if 
$\phi(F^1_1)\subseteq G^1_0=\{0\}$ and 
$\phi(F^2_2)\subseteq G^2_1$ (the other conditions follow from these).

 Let $\phi\in\hoo$ be generic and $S=\ker(\phi)$. It is easy to see that $S$ is of rank $d=1$. Let $\mk=(K^1,K^2)$ be the unique $2$-tuple of subsets of $[2]$ each of cardinality $1$ such that $S\in\cap_{j=1}^2$ $\Omega^o_{K^j}(F^j_{\bull})$. It is easy to see that $K^1=\{1\}$ and $K^2=\{2\}$. By a calculation, we obtain  $\dim(S,V,\mf)=0$. So the filtration in Theorem ~\ref{stonger} is $S\subsetneq V$. Furthermore, the quantity appearing in Inequality
$(\dagger_{\mk}^{\mi})$ is
$\sum_{j=1}^2\sum_{a\in K^j}(n-r+a-i^j_a)-d(n-r)=(2+1)-2=1$
and according to  Theorem ~\ref{stonger}, (iii),  the  dimension of $\hoo$ is $1$ more than its expected dimension (which is $0$). 
\end{example}

The second example is a case where conditions (B) in Theorem ~\ref{MainTe} hold.\begin{example}\label{twoo}
Suppose $s=2$, $r=2$, $n=4$, $I^1=\{1,4\}$ and $I^2=\{2,4\}$. Let $V$ and $Q$ be vector spaces each of rank $2$ and $(\mf,\mg)$ a generic point in $\Fl(V)^2\times\Fl(Q)^2$. A map $\phi$ is in $\hoo$ if  
$\phi(F^1_1)\subseteq G^1_0=\{0\}$ and 
$\phi(F^2_1)\subseteq G^2_1$.
 Suppose $\phi\in\hoo$ is generic and $S=\ker(\phi)$. One can show that (exercise!) $S$ is of rank $d=1$. Let $\mk=(K^1,K^2)$ be the unique $2$-tuple of subsets of $[2]$ each of cardinality $1$ such that $S\in\cap_{j=1}^2$ $\Omega^o_{K^j}(F^j_{\bull})$. It can be shown that  $K^1=\{1\}$ and $K^2=\{2\}$, and by a dimension count,  $\dim(S,V,\mf)=0$. The filtration produced by
Theorem ~\ref{stonger} is $S\subsetneq V$. The quantity appearing in Inequality $(\dagger_{\mk}^{\mi})$ is
$\sum_{j=1}^2\sum_{a\in K^j}(n-r+a-i^j_a)-d(n-r)=(2+0)-2=0$ and   $\hoo$ has rank  equal to its expected dimension (which is $1$). 
\end{example}

The following example produces a filtration with $h=2$. One can similarly create
 examples where $h$ is arbitrarily large.
\begin{example}\label{threee}
Suppose $s=2$, $r=4$, $n=6$, $I^1=\{1,4,5,6\}$ and $I^2=\{2,3,5,6\}$. Let $V$ and $Q$ be vector spaces of ranks $4$ and $2$ respectively  and $(\mf,\mg)$ a generic point in $\Fl(V)^2\times\Fl(Q)^2$. A map $\phi$ is in $\hoo$ if  
$\phi(F^1_1)\subseteq G^1_0=\{0\}$ and 
$\phi(F^2_2)\subseteq G^2_1$.
Pick a generic  $\phi\in\hoo$ and set $S=\ker(\phi)$. It can be shown that  $S$ is of rank $d=2$. Let $\mk=(K^1,K^2)$ be the unique $2$-tuple of subsets of $[2]$ each of cardinality $2$ such that $S\in\cap_{j=1}^2$ $\Omega^o_{K^j}(F^j_{\bull})$. It can be shown  that  $K^1=\{1,4\}$ and $K^2=\{2,4\}$. 

The inductive appeal in the theorem is to the $5$-tuple $(S,V/S,\mf(S),\mf(V/S),\mk)$ which is covered by Example ~\ref{twoo}. Hence the filtration produced by the theorem is of the form
$S^{(2)}\subsetneq S^{(1)}=S\subsetneq V$ ($S^{(2)}\subsetneq S$ is produced by  Example ~\ref{twoo}). By a calculation using Lemma ~\ref{falcon}, it can be shown that 
if   $\mh=(H^1,H^2)$ ($=\mj(2)$ in the theorem)  is the unique $2$-tuple of subsets of $[6]$ each of cardinality $1$ such that $S^{(2)}\in\cap_{j=1}^2$ $\Omega^o_{H^j}(F^j_{\bull})$, then  $H^1=\{1\}$ and $H^2=\{6\}$. It can be checked that $\dim(S^{(2)},V,\mf)=0$. The quantity appearing in Inequality $(\dagger_{\mh}^{\mi})$ is
$\sum_{j=1}^2\sum_{a\in H^j}(n-r+a-i^j_a)-\rk(S^{(2)})(n-r)=(2+0)-2=0.$ The rank of $\hoo$ is  equal to its expected dimension (which is $4$).
\end{example}

The final example is one where $h=1$ and $S^{(1)}=\{0\}$. 
\begin{example}\label{fourr}
Consider the case $s=1$, $r=2$, $n=4$, $I^1=\{2,4\}$. Let $V$ and $Q$ be vector spaces each of rank $2$ and $(\mf,\mg)\in\Fl(V)\times\Fl(Q)$. A map $\phi$ is in $\hoo$ if
$\phi(F^1_1)\subseteq G^1_1$. Suppose $\phi\in\hoo$ is generic. It is easy to see that  $\ker(\phi)=\{0\}$.  The filtration produced by the theorem is $0\subsetneq V$. The unique subset $\mk=K$ of cardinality $0$ of $[2]$ such that $\{0\}\in\Omega^o_K(F^1_{\bull})$ is the empty set and the quantity in Inequality $(\dagger_{\mk}^{\mi})$ is $0$. The rank of $\hoo$ is therefore the expected one ($=2$).
\end{example}
\section{Proof of Theorem ~\ref{MainTe}}\label{FestivusMiracle}
The implication (A)$\Rightarrow$(B) of Theorem ~\ref{MainTe} was shown in Section ~\ref{First}, the implication (B)$\Rightarrow$(A) is proved in  Section ~\ref{second}, and (B)$\Rightarrow$(C) is obvious. To complete the proof of Theorem ~\ref{MainTe}, we need to show (C)$\Rightarrow$(A). For the proof we need to recall the notion of a parabolic vector space from ~\cite{b1}.

 A parabolic vector space $\tilde{V}$ is a $3$-tuple $(V,\mf,w)$, where $V$ is a vector space, $\mf\in\Fl(V)^s$ and $w$ is a function
$$w:\{1,\dots,s\}\times \{1,\dots,\rk(V)\}\to \Bbb{Z}$$ such that if
we let $w^j_l=w(j,l)$, the following holds for each $j=1,\dots,s$:
$$w^j_1\geq w^j_2\geq\dots\geq w^j_{\rk(V)}.$$

Let $S\subseteq  V$ be a non zero subspace of rank $d$.  Let $K^1,\dots,K^s$ be the unique subsets of $[\rk(V)]$ each of cardinality $d$ such that $S\in \cap_{j=1}^s\Omega^o_{K^j}(F^j_{\bull})$. Define the  parabolic slope
$$\mu(S,\tilde{V})=\frac{\sum_{j=1}^s\sum_{a\in K^j} w^j_a}{d}.$$

A parabolic vector space $\tilde{V}$ is said to be semistable if for each subspace $S\subseteq  V$, $\mu(S,\tilde{V})\leq \mu(V,\tilde{V}).$

\subsection{Proof of (C)$\Rightarrow$ (A) in Theorem ~\ref{MainTe}}
Use the notation of Theorem ~\ref{MainTe}. The Chow ring of a Grassmannian does not depend upon the characteristic, therefore we may assume that the base field is of characteristic $0$. We argue by contradiction, assume (C) is true, but (B) is false. Let $\mk=(K^1,\dots,K^s)$, $d$ be as in (B)  and assume Inequality ($\dagger_{\mk}^{\mi}$) fails for this data.

 Let $V$ be an vector space of rank $r$ and $\mf$ a generic point in $\Fl(V)^s$. Define a parabolic vector space $\tilde{V}=(V,\mf,w)$ with  $w^j_a=n-r+a-i^j_a$ for $j=1,\dots,s$ and $a=1,\dots,r$. We have an inequality for the parabolic slope of $V$:
\begin{equation}\label{maltese}
\mu(V,\tilde{V})=\frac{\sum_{j=1}^s\sum_{a=1}^r (n-r+a-i^j_a)}{r}\leq n-r
\end{equation}
because of Inequality ~\ref{conditione} (which is a part of Assumption (C)).

 Since $\prod_{j=1}^s \omega_{K^j}\neq 0$, $\cap_{j=1}^s\Omega^o_{K^j}(F^j_{\bull})$ is non empty. Pick an $S$ in this set. It is easy to see that
$$\mu(S,\tilde{V})=\frac{\sum_{j=1}^s\sum_{a\in K^j} (n-r+a-i^j_a)}{d}.$$
Since Inequality $(\dagger_{\mk}^{\mi})$ fails, and by Inequality ~\ref{maltese}, we have
$$\mu(S,\tilde{V})>n-r\geq\mu(V,\tilde{V}).$$
 The parabolic vector space $\tilde{V}$ is therefore not semistable. Let $T$ be the Harder-Narasimhan maximal contradictor of semistability (see ~\cite{b1}, Lemma 5, proof of Theorem 10). It has the following property:
For any subspace $T'$ of $V$ we have  $\mu(T',\tilde{V})\leq \mu(T,\tilde{V})$ and if $\mu(T',\tilde{V})=\mu(T,\tilde{V})$ then the rank of $T'$ is less than that of $T$.

 Suppose $\tilde{d}=\rk(T)$ and $\ml=(L^1,\dots,L^s)$ the unique $s$-tuple of subsets of $[r]$ each of cardinality $\tilde{d}$ such that $T\in \cap_{j=1}^s\Omega^o_{L^j}(F^j_{\bull})$.

We claim $\cap_{j=1}^s\Omega^o_{L^j}(F^j_{\bull})=\{T\}$: For if  $T'\in\cap_{j=1}^s\Omega^o_{L^j}(F^j_{\bull})$ and $T'\neq T$ then
$\mu(T,\tilde{V})=\mu(T',\tilde{V})$ and the ranks of $T$ and $T'$ are  the same. This contradicts the uniqueness
property of $T$.

Since $\mf$ is generic, by Proposition ~\ref{kleiman} the intersection $\cap_{j=1}^s\Omega^o_{L^j}(F^j_{\bull})$  is transverse at each point of intersection and dense in $\cap_{j=1}^s\Omega_{L^j}(F^j_{\bull})$. So the latter is also $\{T\}$. Hence,
$\prod_{j=1}^s \omega_{L^j}$ is the class of a point in $H^*(\Gr(\tilde{d},r))$ and by  Assumption (C), Inequality $(\dagger_{\ml}^{\mi})$ holds. But this is in contradiction to $\mu(T,\tilde{V})\geq\mu(S,\tilde{V})> n-r$.
\section{Geometric Horn and the  saturation theorem}\label{HornSat}
Recall the definition of the irreducible representations $V_{\lambda}$ of $\GL(r)$ from ~\cite{ful}, Section 3.

For $\lambda =(\lambda_1\geq\dots\geq \lambda_r\geq 0)$ and an integer $l\geq 0$, define $\sigma_{\lambda,l}\in H^{*}(Gr(r,r+l))$ as follows:
\begin{enumerate}
\item $\sigma_{\lambda,l}=0$ if $\lambda_1>l$.
\item If  $\lambda_1\leq l$ define $\sigma_{\lambda,l}$ to be the  cohomology class $\omega_{I(\lambda,l)}$ where $I(\lambda,l)=\{i_1<\dots<i_r\}$ with $i_a=l+a-\lambda_a,$ $a=1,\dots,r$.
\end{enumerate}
Let $R_{+}(\GL(r))$ be the polynomial representation ring of $\GL(r)$ and $l$ a non-negative integer. Define a map of abelian groups $\gamma_l:R_{+}(\GL(r)) \to H^{*}(\Gr(r,r+l))$  by mapping $V_{\lambda}$ to $\sigma_{\lambda,l}$. It is classical that $\gamma_l$ is a homomorphism of rings (see for example ~\cite{young}, Section 9.4).

Let $V_{\lambda}, V_{\mu} ,\dots,V_{\nu}$ be irreducible representations of $\GL(r)$. It is easy to see  that there exist a $V_{\delta}$  with $\delta_1 \leq l$ and an inclusion of $\GL(r)$ representations
$V_{\delta}\subseteq  V_{\lambda}\otimes V_{\mu}\otimes\dots\otimes V_{\nu}$ if and only if $\sigma_{\lambda,l} \cdot\sigma_{\mu,l}\cdot\dots \cdot\sigma_{\nu,l}\neq 0 $ in $H^{*}(\Gr(r,r+l))$.

We now prove Theorem ~\ref{S1}: We need to show that $\sigma_{\lambda,l} \cdot\sigma_{\mu,l}\cdot\dots
\cdot\sigma_{\nu,l}\neq 0 $ in $H^{*}(\Gr(r,r+l))$  if and only if  $\sigma_{N\lambda,Nl}
\cdot\sigma_{N\mu,Nl}\dots\cdot\sigma_{N\nu,Nl}\neq 0 $ in $H^{*}(\Gr(r,r+Nl))$. This is clear from Theorem
~\ref{MainTe} because they are each equivalent to the same set of inequalities (scaled by $N)$. The usual case of
the saturation problem is when the cohomology class of $\sigma_{\lambda,l} \cdot\sigma_{\mu,l}
\cdot\dots\cdot\sigma_{\nu,l}$ is of the top degree in $H^{*}(Gr(r,r+l))$ (namely $(2)rl$).
\section{Genericity}\label{whatdreamsmaycome}
\subsection{Stratification and Universal families}\label{hamlett}
Let $\mi=(I^1,\dots,I^s)$, $V$ and $Q$ be as in Definition~\ref{deftwo}. Let
$\mk=(K^1,\dots,K^s)$ be a $s$-tuple of subsets of $[r]$ each of cardinality $d$. We consider the following ``universal objects'' (the scheme structures on these sets are from the theory of
determinantal schemes and will be discussed in Appendix ~\ref{determinant}).
\begin{enumerate}
\item[(A)] Define  $\mathfrak{H}_{\mi}(V,Q,\mk)$  to be the scheme
over  $\Fl(V)^s\times \Fl(Q)^s$ whose fiber over $(\mf,\mg)$ is
(same as the space given by Equation ~\ref{latenight})
\begin{equation}\label{goo}
\{\phi\in \hoo\mid \rk(\ker(\phi))=d,\text{ }\ker(\phi)\in \cap_{j=1}^s\Omega^o_{K^j}(F^j_{\bull})\}.
\end{equation}
\item[(B)] Define $\muu_{\mk}(V)$ to be the scheme over $\Fl(V)^s$ whose fiber over $\mf\in \Fl(V)^s$ is  $\cap_{j=1}^s\Omega^o_{K^j}(F^j_{\bull})\subseteq \Gr(d,V)$.
\end{enumerate}

The scheme structures are such that the induced flag constructions on the subbundle and the quotient bundle are algebraic. For example,
in $\muu_{\mk}(V)$, the universal subbundle $\ms$ has $s$ (induced) complete filtrations by subbundles (and the same for the universal quotient). The following proposition is proved in the appendix.
\begin{proposition}\label{forlater}
\begin{enumerate}
\item $\mhh_{\mi}(V,Q,\mathcal{K})$  and $\mukv$ are smooth  and irreducible schemes.
\item If $\ham\neq\emptyset$, the natural morphism $p:\ham\to \mukv$ which maps $(\phi,\mf,\mg)$ to $(\ker(\phi),\mf)$, is smooth and surjective.
\end{enumerate}
\end{proposition}
We give the essential ideas behind Proposition ~\ref{forlater}, postponing the complete proofs to the appendix. To
make the arguments below rigorous, it is necessary to use Grothendieck's
language of functors and representability (which we do in the appendix).

Let $\ms$ be the
universal subbundle of $V\tensor \mathcal{O}_{\Gr(d,V)}$ on
$\Gr(d,V)$ and $\Fl(\ms)\to \Gr(d,V)$ the complete flag bundle of
$\ms$. Consider the $s$-fold fiber product
$$\Fl(\ms)^s=\Fl(\ms)\times_{\Gr(d,V)}\dots\times_{\Gr(d,V)}\Fl(\ms)$$
 Points of $\Fl(\ms)^s$ are pairs $(S,\mathcal{R})$ where
 $S\in \Gr(d,V)$ and $\mathcal{R}\in\Fl(S)^s$.
 Consider the morphism $\mukv\to\Fl(\ms)^s$ obtained by taking the induced
flags on the subspace. The fiber of this morphism over $(S,\mathcal{R})\in\Fl(\ms)^s$ is the set of choices of $\mf\in \Fl(V)^s$ so
that $S\in\cap_{j=1}^s\Omega^o_{K^j}(F^j_{\bull})$ and the collection of flags induced on $S$ by $\mf$ is
$\mathcal{R}$. Each coordinate of $\mf$ can be selected ``separately and independently''. The set of choices of each coordinate
is  parameterized by an open subset of the top most member of  a tower of projective bundles over $\spec{\kappa}$ (for each $j$, consider
Lemma ~\ref{basic11} with $\pl=V$,  ${\pl}^l=R^j_l$,  $a_l=k^j_l$, $l=1,\dots,d$, $k=d$, $\pl^{d+1}=V$ and $\pl^d =S$). It
follows that $\mukv\to\Fl(\ms)^s$ is a ``locally trivial fiber bundle'' with a smooth fiber. This shows the
smoothness and irreducibility of $\mukv$.

Let $\mathcal{A}=\mukv$. Denote by $\ms_{\mathcal{A}}$ the subbundle of $V_{\mathcal{A}}=V\tensor\mathcal{O}_{\ma}$ on $\mathcal{A}$ induced from $\ms$. Form the total space $Y$ of the vector bundle $\shom(V_{\ma}/\ms_{\ma},Q_{\ma})$ over $\ma$ and consider the open part $Y^o$ of it when the universal morphism $V_{Y}/\ms_{Y}\to Q_{Y}$ is an injection (on the fibers). Points of $Y^o$ are pairs $(\mf,\phi)$ where
$\mf\in\Fl(V)^s$, $\phi\in \home(V,Q)$ such that
$\rk(\ker(\phi))=d$ and  $\ker(\phi)\in \cap_{j=1}^s\Omega^o_{K^j}(F^j_{\bull})$. It is easy to see that the natural map $\ham \to Y^o$ is a ``locally trivial fiber bundle'' with smooth fibers: The fiber over $(\mf,\phi)$ corresponds to the set of choices of $\mg\in\Fl(Q)^s$  so that $\phi\in\hoo$.
Each coordinate of $\mg$ can be selected ``separately and independently'': The defining condition for $\phi\in \hoo$ is
$$G^j_{i^j_l-l}\supseteq \phi(F^j_l)$$
for $l=1,\dots,r$ and $j=1,\dots,s$. The rank of $\phi(F^j_l)$ is $l-t$ if $k^j_t\leq l<k^j_{t+1}$. For each $j$, apply Lemma ~\ref{basic11} with $\pl=Q$, ${\pl}^l=\phi(F^j_l)$ and $a_l=i^j_l-l$,  $l=1,\dots,r$. Proposition ~\ref{forlater} follows from these considerations (the map $Y^o\to \mathcal{A}$ is also a fiber bundle with smooth fibers).

The following elementary result is a part of Lemma ~\ref{basic1} (for $Z=\spec{\kappa}$).
\begin{lemma}\label{basic11}
Let $\pl$ be a vector space of rank $f$. Assume that $\pl$ is filtered by a series of vector subspaces:
$$0={\pl}^{0}\subseteq {\pl}^{1}\subseteq  \dots \subseteq  {\pl}^k\subseteq {\pl}^{k+1}={\pl}.$$
 Let $0=a_0\leq a_1\leq\dots\leq a_k\leq f$ be nonnegative integers. Consider the variety $A$ of  complete filtrations by subspaces:
$$0\subsetneq F_1\subsetneq\dots\subsetneq F_f={\pl}$$
so that
$$
F_{a_l}\supseteq {\pl}^{l}\text{ for }l=1,\dots,k
$$
 Then, $A$ is the topmost element in  a tower of projective bundles over $\spec{\kappa}$. That is, there exists a sequence of morphisms
$$A=A_{h+1}\leto{p_h} A_{h}\leto{p_{h-1}}A_{h-1}\leto{p_{h-2}}\dots\leto{p_0}A_0=\spec{\kappa}$$
 such that for $u=1,\dots,h$, $A_{u+1}$ is the projective bundle associated to a vector bundle $\mw_u$ on $A_u$ and $p_u$ the projection map.

If we impose the further condition:
$$\rk(F_{b}\cap {\pl}^{k})=l$$ for $a_l\leq b<a_{l+1}$
and $l=1,\dots,k$, then the corresponding space of filtrations is an open (possibly empty) subset of $A$.
\end{lemma}
\subsubsection{Genericity of induced structures, the main idea:}
We give the essential ideas behind the genericity statements of Section ~\ref{strategem2}, postponing the complete proofs to Section ~\ref{dreams}.

Genericity statements usually follow from the properties of the relevant ``universal families''. For example, using Proposition ~\ref{forlater}, one deduces that if $U$ is a nonempty open subset (points satisfying some desired ``open'' property for instance) of $\mukv$, there is a nonempty open subset of points $\tilde{U}$ of $\ham$ which have for their image in $\mukv$, a point in $U$. By Lemma ~\ref{nice}, a general fiber of $\ham\to\Fl(V)^s\times\Fl(Q)^s$ has a dense intersection with $\tilde{U}$. We can therefore say (somewhat imprecisely) that if $(\mf,\mg)$ is a generic point of $\Fl(V)^s\times\Fl(Q)^s$ and $\phi$ a generic element of the fiber of $\ham\to \Fl(V)^s\times \Fl(Q)^s$ over $(\mf,\mg)$ (given by Equation $~\ref{goo}$), then
the induced point $(\ker(\phi),\mf(\ker(\phi)))\in\mukv$ is generic.

Suppose that we are given a $\GL(\kappa^d)\times\GL(\kappa^{r-d})$ invariant nonempty open subset  $\mathfrak{O}$ of $\Fl(\kappa^d)^s\times\Fl(\kappa^{r-d})^s$.
An example of the nonempty open property  of  points $(S,\mf)$ in $\mukv$ could then be that the induced pair of flags $(\mf(S),\mf(V/S))$ is carried to a point of $\mathfrak{O}$ upon some (hence any, because of the invariance property of $\mathfrak{O}$) identification $S\to \kappa^d$, $V/S\to \kappa^{r-d}$. That the corresponding open subset $U$ of $\mukv$ is nonempty follows from the assumed nonemptiness of $\mathfrak{O}$: Take a $T$ in $\Gr(d,V)$ and $\mh\in \Fl(T)^s$ and  $\mrr\in \Fl(V/T)^s$, such that  $(\mh,\mrr)$ has the required property. By Lemma ~\ref{flags}, we can find a $\mf\in \Fl(V)^s$ such that $(T,\mf)\in\mukv$ and the ordered collection of flags induced by it on $T$ and $V/T$ are $\mh$ and $\mrr$ respectively.

 The above discussion  connects with the genericity property (G3) of Section ~\ref{strategem2}. We make these ideas precise in the next section.

 The following statement on the genericity of induced flags is proved by the same methods as above (for the proof use the irreducibility of $\muiw$).  We are not going to use it in this paper. 
\begin{proposition}\label{neww}
Let  $W$ be a vector space of rank $n$ and $\mi=(I^1,\dots,I^s)$ a $s$-tuple of  subsets of $[n]$ each of cardinality $r$. For generic $\me\in\Fl(W)^s$, $\cap_{j=1}^s\Omega^o_{I^j}(E^j_{\bull})$ has an open dense subset $O$ so that for $V\in O$, $\me(V)$ is generic for intersection theory in $V$ (that is, $\me(V)\in B(V)$, see Section ~\ref{strategem2} for the definition of $B(V)$).
\end{proposition}

\begin{lemma}\label{nice} Let $X,$ $Y$ be irreducible algebraic varieties, $f:X\to Y$ a morphism and $U(X)\subseteq  X$ a non-empty open subset. Then, there exists a nonempty open subset $U(Y)$ of $Y$ such that
\begin{enumerate}
\item Either $f^{-1}(U(Y))=\emptyset$ or $f:f^{-1}(U(Y))\to U(Y)$ is flat and surjective.
\item For $y\in U(Y)$, $f^{-1}(y)\cap U(X)$ is
dense in $f^{-1}(y)$ (which could be empty).
\end{enumerate}
\end{lemma}
\begin{proof}
 Let $Z=X-U(X)$.  Each irreducible component of $Z$ is of dimension strictly less than the dimension of $X$. Using generic flatness (see ~\cite{eisen}, Theorem 14.4), find a non-empty  open $U\subseteq  Y$ such that
\begin{enumerate}
\item[(a)]$f:f^{-1}(U)\to U$ is flat.
\item[(b)]$f\mid_Z:f^{-1}(U)\cap Z\to U$ is flat.
\end{enumerate}
For $y\in U$ because of the flatness conditions above, each irreducible component of the fiber $f^{-1}(y)$ is of dimension $\dim(X)-\dim(Y)$ and each irreducible component of $f^{-1}(y)\cap Z$ has dimension strictly less than $\dim(X)-\dim(Y)$. Therefore, if $y\in U$, $f^{-1}(y)\cap U(X)$ is
dense in $f^{-1}(y)$ (which could be empty).

If $f^{-1}(U)=\emptyset$ set $U(Y)=U$ and if  $f^{-1}(U)\neq \emptyset$, set $U(Y)=f(f^{-1}(U))$ which is then a nonempty open set of $U$ (flat maps between varieties are open, see ~\cite{hart}, III, Exercise 9.1). It is easy to see that $U(Y)$ satisfies the required properties.
\end{proof}
\subsection{Proof of the genericity statement from Section ~\ref{strategem2}}\label{dreams}
Let us recall some open subsets in the parameter spaces considered before:

Recall from Section ~\ref{strategem2}, the nonempty Zariski open subset $B(V)\subseteq  \Fl(V)^s$ for every nonzero vector space $V$. It is easy to see that if $V\to V'$ is an isomorphism, $B(V)$ maps bijectively to $B(V')$ under the natural bijection $\Fl(V)^s\to\Fl(V')^s$ (``functoriality for isomorphisms'').

 Recall the open subset $O(V,Q,\mi)\subseteq  \Fl(V)^s\times \Fl(Q)^s$ from  Section ~\ref{strategem2}. Let $\mk=\mk(\mi)$. We note the following properties:
\begin{enumerate}
\item The map $\ham\to \Fl(V)^s\times\Fl(Q)^s$ is flat over $O(V,Q,\mi)$ because  over $O(V,Q,\mi)$ it is surjective, has constant fiber dimension and the source and target are both smooth (see ~\cite{mat}, Theorem 23.1).
\item $O(V,Q,\mi)$ is  ``functorial'' for isomorphisms: Given isomorphisms $V\to V'$ and $Q\to Q'$, the induced map $\Fl(V)^s\times\Fl(Q)^s\to \Fl(V')^s\times\Fl(Q')^s$ takes $O(V,Q,\mi)$ bijectively to $O(V',Q',\mi)$ (This follows from Kleiman's proof ~\cite{kl} of Proposition ~\ref{kleiman}).
\end{enumerate}
The reader unwilling to look at Kleiman's method of proof of Proposition ~\ref{kleiman}, may just redefine $O(V,Q,\mi)$ to be the largest open subset of $\Fl(V)^s\times\Fl(Q)^s$ such that given any $s$-tuple $\mj$ of subsets of $[r]$ each of the same cardinality,  the following property is satisfied: Denoting by $\pi_{\mj}$ the natural map $\mhh_{\mi}(V,Q,\mj)\to \Fl(V)^s\times\Fl(Q)^s$,
either $\pi_{\mj}^{-1}(O(V,Q,\mi))=\emptyset$ or $\pi_{\mj}:{\pi}_{\mj}^{-1}(O(V,Q,\mi))\to\Fl(V)^s\times\Fl(Q)^s$ is flat and surjective (see Lemma ~\ref{nice} to show the nonemptiness of $O(V,Q,\mi)$). Functoriality for isomorphisms follows from this construction.

We now construct $A(V,Q)$ (whose existence was claimed in Section ~\ref{strategem2}, satisfying properties (G1), (G2) and  (G3)) which is ``functorial'' for isomorphisms (and hence $A(V,Q)$ is $\GL(V)\times \GL(Q)$ invariant). We do this inductively.

Notice that the conditions (G2) and (G3) are requirements over a finite set of possible $s$-tuples $\mi=(I^1,\dots,I^s)$ (given $V$ and $Q$). Assume that we have constructed  $A(V,Q)$ for all $V$ of rank $<r_0$ and arbitrary rank of $Q$.  We now construct $A(V,Q)$ if the rank of $V$ is $r=r_0$. The base case of induction is by setting $r_0=1$ in the argument below.  The requirement on $A(V,Q)$ is over all possible choices of $\mi$. We will construct a nonempty open subset $A(V,Q,\mi)\subseteq  \Fl(V)^s\times \Fl(Q)^s$, which satisfies (G2) and (G3) for $\mi$, and the requirement (G1).

 We may then finally define $A(V,Q)=\cap_{\mi}A(V,Q,\mi)$ (a finite intersection). It is immediate that $A(V,Q)$ satisfies (G1), (G2) and (G3). The functorial property for isomorphisms will be  clear from the construction of $A(V,Q,\mi)$.
\subsubsection{Construction of $A(V,Q,\mi)$}\label{openity}
Let $\mi$ be a $s$-tuple as above. Let $d=d(\mi)$ and $\mk=\mk(\mi)=(K^1,\dots,K^s)$. If $d=0$ or $d=r$: Then just take $A(V,Q)=O(V,Q,\mi)\cap [B(V)\times B(Q)]$. Now assume $0<d<r$ and consider the morphisms
$$\xymatrix{ &   \ham\ar[dl]^{p}\ar[dr]^{\pi}\\
         \mukv  &  &     \Fl(V)^s\times \Fl(Q)^s }$$
Let $U_{o}$ be the subset of $\mukv$ formed by points $(T,\mf)$  such that $(\mf(T),\mf(V/T))\in A(T,V/T)$.  We first show that $U_{o}$ is a nonempty  open
subset of $\mukv$.
\begin{enumerate}
\item {\bf Non-emptiness:} By induction, for any $d$-dimensional subspace $T\subseteq  V$, we can find $\mh\in \Fl(T)^s$ and  $\mrr\in \Fl(V/T)^s$, such that  $(\mh,\mrr)\in A(T,V/T)$. By Lemma ~\ref{flags}, we can find a $\mf\in \Fl(V)^s$ such that $(T,\mf)\in\mukv$ and the ordered collection of flags induced by it on $T$ and $V/T$ are $\mh$ and $\mrr$ respectively.  Hence $U_{o}$ is non empty.
\item {\bf Openness:} For every $P\in\mukv$, there is a neighborhood $U(P)$ of $P$ in $\mukv$ and a morphism
(non-canonical) $\gamma:U(P)\to \Fl(E)^s\times \Fl(M)^s$ where $E$ and $M$ are vector spaces of ranks $d$ and $r-d$ respectively, such that $E$, $M$ (and their filtrations) pull back to the universal subbundle and the quotient bundle on $\mukv$ respectively (and their induced filtrations). It is easy to see that $U_{o}\cap U(P) =\gamma^{-1}(A(E,M))$ and hence $U_{o}$ is an open
subset of $\mukv$.
\end{enumerate}
 Let $\breve{U}_{o}=p^{-1}(U_{o})$. By Proposition ~\ref{forlater}, $p$ is surjective, therefore ${\breve{U}}_{o}\neq \emptyset$. Since $\ham$ is irreducible, $\pi^{-1}(O(V,Q,\mi))\cap {\breve{U}}_{o}$ is a non-empty open set. We now define
 $$A(V,Q,\mi)=\pi({\breve{U}}_{o})\cap O(V,Q,\mi)\cap (B(V)\times B(Q))$$
 and is a non-empty open subset of $\Fl(V)^s\times \Fl(Q)^s$ (recall that
  $\pi$ is flat over $O(V,Q,\mi)$ and that flat maps are open).

It is easy to see that $A(V,Q,\mi)$
 fulfils the requirements (G2) and (G3) for $\mi$,  and (G1):
\begin{itemize}
\item (G1) is satisfied since $A(V,Q,\mi)\subseteq  B(V)\times B(Q)$.
\item (G2) is seen to be satisfied for $\mi$ from the inclusion
$A(V,Q,\mi)\subseteq  O(V,Q,\mi)$.
\item (G3) is satisfied for $\mi$ because if $(\mf,\mg)\in A(V,Q,\mi)$, then $\breve{U}_{o}\cap \pi^{-1}(\mf,\mg)$ is a nonempty open subset of $\pi^{-1}(\mf,\mg)$. Since $(\mf,\mg)\in O(V,Q,\mi)$,  $\pi^{-1}(\mf,\mg)$ is a dense open subset of the vector space $\hoo$ and is therefore  irreducible. It follows that  $\breve{U}_{o}\cap \pi^{-1}(\mf,\mg)$ is a nonempty dense open subset of $\pi^{-1}(\mf,\mg)$ and hence (G3) holds.
\end{itemize}
\section{Transversality in Schubert Calculus}
\begin{theorem}\label{transversality}
Let  $W$ be a vector space of rank $n$ and $\mi=(I^1,\dots,I^s)$ a $s$-tuple of  subsets of $[n]$ each of cardinality $r$. Then, for generic $\me\in\Fl(W)^s$, $\cap_{j=1}^s\Omega^o_{I^j}(E^j_{\bull})$ has an open dense subset of points  where the intersection is transverse.
\end{theorem}
\begin{proof}
Let $q:\muiw\to \Fl(W)^s$. If $q$ is not dominant, then there is nothing to prove. Assume  $q$ is dominant, then clearly consideration of a general fiber of $q$ leads us to the conclusion $\prod_{j=1}^s  \omega_{I^j}\neq 0$ (See Lemma ~\ref{k1}).

By Proposition ~\ref{forlater}, $\muiw$ is smooth and irreducible. Consider the   subset $U$ of $\muiw$ formed by  points $(T,\me)$ satisfying the following property: $\Omega^o_{I^1}(E^1_{\bull}),\dots,\Omega^o_{I^s}(E^s_{\bull})$ meet transversally at $T$.

We claim that  $U\subseteq \muiw$ is a nonempty open subset. The openness follows from the observation that $U$ is the smooth locus of $q$. Alternately, as in Section ~\ref{openity}, for every $P\in\muiw$, there is a neighborhood $U(P)$ of $P$ in $\muiw$ and a morphism
(non-canonical) $\gamma:U(P)\to \Fl(E)^s\times \Fl(M)^s$ where $E$ and $M$ are vector spaces of ranks $r$ and $n-r$ respectively, such that $E$, $M$ (and their filtrations) pull back to the universal subbundle and the quotient bundle on $\muiw$ respectively (and their induced filtrations). Let $\tilde{U}$ be the open subset of $\Fl(E)^s\times\Fl(M)^s$ (see Lemma ~\ref{april1}, (2)) formed by points $(\mf,\mg)$ such that the rank of $\home_{\mi}(E,M,\mf,\mg)$ is $\dim(\Gr(r,n))-\sum_{j=1}^s \codim(\omega_{I^j})$. It is easy to see that $U\cap U(P) =\gamma^{-1}(\tilde{U})$ and hence $U$ is an open subset of $\muiw$.

We now show that $U$ is nonempty:  Since $\prod_{j=1}^s  \omega_{I^j}\neq 0$, Property (A) holds in Theorem ~\ref{MainTe}. By Section ~\ref{First}, Property (B) holds in Theorem ~\ref{MainTe}. By Section ~\ref{second}, the condition $(\gamma)$ of Proposition ~\ref{april2} holds. By Proposition ~\ref{april2}, $(\gamma)\Rightarrow (\beta)$, and we find a point of $\muiw$ with the desired property.

We now make an appeal to Lemma ~\ref{nice},  with $X=\muiw$, $Y=\Fl(W)^s$ and $U(X)=U$ and conclude the proof.
\end{proof}
The transversality part of this theorem is weaker in characteristic $0$ than the statement obtained by using
Kleiman's transversality theorem (Proposition  ~\ref{kleiman}). In characteristic $0$ for any $s$-tuple $\mi$ as above, and $\me\in \Fl(W)^s$ a generic point, $\Omega^o_{I^1}(E^1_{\bull})$, $\dots$, $\Omega^o_{I^s}(E^s_{\bull})$  meet transversally at each point of intersection. However for enumerative problems where the expected dimension is $0$, Theorem ~\ref{transversality}
is just as strong:
\begin{corollary}\label{sottilee}
Let  $W$ be a vector space of rank $n$. Consider a  $s$-tuple $\mi=(I^1,\dots,I^s)$ of subsets of $[n]$ each of cardinality $r$, such that
$$\dim(\Gr(r,n))-\sum_{j=1}^s\codim(\omega_{I^j})=0$$
Then, for generic $\me\in\Fl(W)^s$,
 $\cap_{j=1}^s \Omega_{I^j}(E^j_{\bull})$, is a reduced zero dimensional scheme (possibly empty).
\end{corollary}
\begin{proof}
By Proposition ~\ref{kleiman}, for  generic $\me\in\Fl(W)^s$, $\cap_{j=1}^s \Omega_{I^j}(E^j_{\bull})$ is zero dimensional and $\cap_{j=1}^s \Omega^o_{I^j}(E^j_{\bull})$ is dense in it. That is, $\cap_{j=1}^s \Omega_{I^j}(E^j_{\bull})=\cap_{j=1}^s \Omega^o_{I^j}(E^j_{\bull})$. We now apply Theorem ~\ref{transversality} and conclude the proof.
\end{proof}
 \appendix\section{Determinantal Schemes}\label{determinant}
\subsection{Standard Theory}
Let $X$ be a scheme and let $V$ and $Q$ be free $\mathcal{O}_X$-modules of ranks $e$ and $f$, respectively. Let $\phi:V\to Q$ be a morphism. The rank of $\phi$ is said to be $r$ if the cokernel is a locally free sheaf of rank $f-r$. It is easy to see that if the cokernel of $\phi$ is locally free, then so are the image and the kernel.

With $X$, $V$, $Q$ and  $\phi$ as above, the $r$th degeneracy locus $D_r(\phi)$ of $\phi$ is the closed subset of all points $x\in X$ such that $\phi\tensor k(x)$ has rank less than or equal to $r$. We put a structure of a closed subscheme on $D_r(\phi)$ by writing it as the scheme of zeroes of the morphism
$$\bigwedge^{r+1}V\to \bigwedge^{r+1}Q.$$
Notice that the subschemes $Z_r(\phi)=D_r(\phi)-D_{r-1}(\phi)$ partition $X$ into a disjoint union of locally closed subschemes.

\begin{lemma} With notation as above, if $T=Z_r(\phi)$, the map $\pi_T:V_T\to Q_T$, has rank $r$. That is, the cokernel is locally free of dimension $f-r$.
\end{lemma}
\begin{proof}See ~\cite{eisen}, Corollary 20.5, Proposition 20.8.
\end{proof}
Given $X$, $V$, and $Q$ as above we form the scheme $\home(V,Q)$ over $X$. Write $\pi:\home(V,Q)\to X$ for the structure map. There is a natural map $\phi:\pi^{*}V\to \pi^{*}Q$. we write $\home_r(V,Q)$ for the locally closed subscheme $Z_r(\phi)$. This represents a functor:
\begin{lemma}
Let $X$ be a Noetherian scheme, then $\home_r(V,Q)$ represents the functor
$$T\leadsto \{\psi:V_T\to Q_T\mid\operatorname{rank}(\psi)=r\}.$$ 
\end{lemma}
\begin{proof}
This is standard, for example see \cite{bro}, Prop 2.1.
\end{proof}
\subsection{Universal families}
We will use Grothendieck's language of functors and representability to discuss Grassmannians and flag varieties (see ~\cite{EGA}, \S 9). All schemes considered in this section are over $\text{Spec}(\kappa)$.
\begin{lemma}\label{basic1}
Let ${\tpl}$ be a vector bundle of rank $f$ on a scheme $Z$. Assume that ${\tpl}$ is filtered by a series of subbundles:
$$0={\tpl}^{0}\subseteq {\tpl}^{1}\subseteq  \dots \subseteq  {\tpl}^{k+1}={\tpl}.$$
 Let $0=a_0\leq a_1\leq\dots\leq a_k\leq f$ be nonnegative integers. Consider the functor $A$ whose value over a  scheme $T$ over $Z$ is the set  of  complete filtrations by subbundles:
$$0\subsetneq \mathfrak{F}_1\subsetneq\dots\subsetneq \mathfrak{F}_f={\tpl}_T$$
so that
\begin{equation}\label{starre}
\mathfrak{F}_{a_l}\supseteq {\tpl}^{l}_T\text{ for }l=1,\dots,k
\end{equation}
 This functor is representable by the topmost element in  a tower of  Grassmann bundles (actually bundles of projective spaces) over $Z$.

Furthermore, consider the functor $\tilde{A}$ which is the same as the functor $A$ but for the additional condition that
the morphism
$$\tpl^{k}\to \tpl/\mathfrak{F}_b$$
 is of rank $\rk(\tpl^{k})-l$ for $a_l\leq b<a_{l+1}$ and  $l=1,\dots,k$. Then, the functor $\tilde{A}$ is representable by an open subset of the representing scheme of $A$.
\end{lemma}
\begin{proof} The result is best understood at the ``level of points''.  Assume $Z=\spec(\kappa)$. We start by choosing ${\mathfrak{F}}_{f-1}$ so that it contains ${\tpl}^{t}$ whenever $a_t\leq f-1$. That is, we choose ${\mathfrak{F}}_{f-1}$, a subspace of ${\tpl}$ which contains ${\tpl}^{t_1}$ where $t_1$ is the largest $t$ satisfying $a_t\leq f-1$. There is a projective space worth of possible choices.  Having chosen ${\mathfrak{F}}_{r-1}$, we choose ${\mathfrak{F}}_{f-2}$ contained in ${\mathfrak{F}}_{f-1}$ and containing ${\tpl}^{t_2}$
 where $t_2$ is the largest $t$ satisfying $a_t\leq f-2$ (the choices are again parameterized by points in a projective space) and so on. The condition (for the functor $\tilde{A}$): $\rk(\mathfrak{F}_b\cap{\tpl}^{k+1})=l$ for $a_l\leq b<a_{l+1}$ and  $l=1,\dots,k$ clearly imposes an open condition. We cast this in the language of functors:

For $j=1,\dots,f$ let
$$c_j=\max\{l\mid a_l\leq j\}.$$
Condition (~\ref{starre}) can be rewritten as
$${\tpl}^{c_j}_T\subseteq  {\mathfrak{F}}_j,\text {for }j=1,\dots,f.$$
 Let  $X_{f-1}=\Gr(f-1-\rk({\tpl}^{c_{f-1}}),\frac{{\tpl}}{{\tpl}^{c_{f-1}}})$ and $p_{f-1}:X_{f-1}\to Z$. On $X_{f-1}$ there is a natural bundle $\mt_{f-1}$ of rank $f-1$ such that there are natural inclusions of subbundles
$$p_{f-1}^*({\tpl}^{c_{f-1}})\subseteq \mt_{f-1}\subseteq p_{f-1}^*{\tpl}.$$
Let $p_{f-2}:X_{f-2}\to X_{f-1}$ be the Grassmann bundle (these are bundles of projective spaces, or identity morphisms)
$$\Gr(f-2-\rk({\tpl}^{c_{f-2}}),\frac{\mt_{f-1}}{p_{f-1}^*({\tpl}^{c_{f-2}})}  )$$
and so on obtaining a tower of Grassmann bundles 
$$X_1\leto{p_1}X_2\leto{p_2}\dots\leto{p_{f-2}}X_{f-1}\leto{p_{f-1}} Z.$$
Let $q:X_1\to Z$ be the composition $p_{f-1}\circ p_{f-2}\circ\dots\circ p_1$. It is easy to see that $X_1$ represents the functor $A$ ($q^*(\tpl)$ has the required  universal filtration coming from the pullback of $\mt_{\bull}$). The conclusion for the
functor $\tilde{A}$ is immediate.
\end{proof}
We use the notation of Section ~\ref{hamlett}.
For a scheme $T$, let $B(\mk,V)[T]$ be the set  of ordered data consisting of:
\begin{enumerate}
\item[(a)] A subbundle $\ms\subseteq  \mv_T=V\tensor_{\kappa} \mathcal{O}_T$ of rank $d$.
\item[(b)] The data of $s$ complete filtrations of $\mv_T$. That is, for each $1\leq j\leq s$, a filtration by subbundles:
$$0\subsetneq F^j_1(\mv_T)\subsetneq\dots\subsetneq F^j_r(\mv_T)=\mv_T$$
with $\rk(F^j_a(\mv_T))=a$ for $a=1,\dots,r$.
\end{enumerate}
This data is required to satisfy the condition: The rank of $\ms\to \mv_T/F^j_a(\mv_T)$ is $d-t$ for $k^j_t\leq a<k^j_{t+1}$, $t=0,\dots,d$, $j=1,\dots,s$.

In the above situation, the bundle $\ms$ gets $s$
induced complete filtrations by subbundles as follows: Let ${\ram}^j_t=\ker(\ms\to \mv/F^j_{k^j_t}(\mv_T))$ for $t=1,\dots,d$ and $j=1,\dots,s$. So we see that: The functor $B(\mk,V)$ is naturally isomorphic to the functor $\tilde{B}(\mk,V)$, which assigns to every scheme $T$, the set of ordered data consisting of
\begin{enumerate}
\item[(a)] A subbundle $\ms\subseteq  \mv_T = V\tensor_{\kappa} \mathcal{O}_T$ of rank $d$.
\item[(b)] For each $1\leq j\leq s$, a complete filtration by subbundles:
$$0\subsetneq F^j_1(\mv_T)\subsetneq\dots\subsetneq F^j_r(\mv_T)=\mv_T.$$
\item [(c)] For each $1\leq j\leq s$, a complete filtration by subbundles:
 $$0\subsetneq {\ram}^j_1\subsetneq\dots\subsetneq {\ram}^j_d=\ms.$$
\end{enumerate}
This data is required to satisfy
\begin{enumerate}
\item[($T_1$)] ${\ram}^j_t\subseteq  F^j_{k^j_t}(\mv_T)$ for $j=1,\dots,s$ and $t=1,\dots,d$.
\item[($T_2$)] The rank of $\ms\to \mv/F^j_a(\mv_T)$ is $d-t$ for $k^j_t\leq a<k^j_{t+1}$, $t=0,\dots,d$, $j=1,\dots,s$.
\end{enumerate}
 For a scheme $X$ and a vector bundle $\mt$ over it, let $\Fl(\mt)\to X$ denote the scheme of complete flags on the fibers of $\mt$.
\begin{lemma}\label{hamlet}Let $\ms$ be the universal subbundle of $V\tensor\mathcal{O}_{\Gr(d,V)}$ on $\Gr(d,V)$.
The functor $B(\mk,V)$ is represented by an open subset of the topmost element
 of a tower of Grassmann bundles over
$$X=\Fl(\ms)\times_{\Gr(d,V)}\times\dots\times_{\Gr(d,V)}\Fl(\ms)$$
with $\Fl(\ms)$ repeating $s$ times. The representing scheme is therefore smooth and irreducible. At the level of sets the representing scheme  coincides with the set $\mukv$ from Section ~\ref{whatdreamsmaycome}.
\end{lemma}
\begin{proof}
It is easy to see that there is a natural transformation  from $\tilde{B}(\mk,V)$($\leto{\sim}B(\mk,V)$) to the functor of points of $X$. The condition ($T_1$) and Lemma ~\ref{basic1} give a natural transformation from $\tilde{B}$ to the (functor of points of the)  topmost element $Y$ of a tower of Grassmann bundles over $X$ (``we have to chose the flags on $V$''). This scheme $Y$ carries a universal family for the functor $\tilde{B}(\mk,V)$ except that the condition ($T_2$) may not be valid for this universal family.  We can write the condition $(T_2)$ as
\begin{enumerate}
\item[($T_2'$)]The morphism $\ms/{\ram}^j_t\to \mv/F^j_a(\mv_T)$ is of maximal rank ($=d-t$) for $k^j_t\leq a<k^j_{t+1}$, $t=0,\dots,d$, $j=1,\dots,s$.
\end{enumerate}
It is easy to see that ($T_2'$) is an open condition on this family $Y$. Therefore $\tilde{B}(\mk,V)$ is represented by an open subset of $Y$. The assertion on the points of $\mukv$ is immediate from the definitions.
\end{proof}
 For a scheme $T$, let $\mv_T=V\tensor_{\kappa} \mathcal{O}_T$ and $\mq_T=Q\tensor_{\kappa} \mathcal{O}_T$. Define a functor $H_{\mi}(V,Q,\mk)$ which assigns to a scheme $T\text{ over Spec}(\kappa)$ the set whose elements are the following (ordered) data (conditions (a) and (b) below are the same as the conditions in the functor $B(\mk,V)$):
\begin{enumerate}
\item[(a)] A subbundle $\ms\subseteq  \mv_T = V\tensor_{\kappa} \mathcal{O}_T$ of rank $d$.
\item[(b)] For each $1\leq j\leq s$, a complete filtration by subbundles:
$$0\subsetneq F^j_1(\mv_T)\subsetneq\dots\subsetneq F^j_r(\mv_T)=\mv_T $$
 satisfying the following condition: The rank of $\ms\to \mv_T/F^j_a(\mv_T)$ is $d-t$ for $k^j_t\leq a<k^j_{t+1}$, $t=0,\dots,d$, $j=1,\dots,s$.
\item [(c)]For each $1\leq j\leq s$, a complete filtration by subbundles:
$$0\subsetneq G^j_1(\mq_T)\subsetneq\dots\subsetneq G^j_{n-r}(\mq_T)=\mq_T.$$
\item[(d)] A morphism $\phi:\mv_T\to \mq_T$ of rank $r-d$ and kernel $\ms$ such that for $j=1,\dots,s$ and $a=1,\dots,r$ we have
$$\phi(F^j_a(\mv_T)) \subseteq  G^j_{i^j_a-a}(\mq_T).$$
\end{enumerate}
\begin{lemma}
\begin{enumerate}
\item The functor $H_{\mi}(V,Q,\mk)$ is representable by a smooth scheme. This representing scheme coincides  at the level of points with the set $\mhh_{\mi}(V,Q,\mk)$ in Section ~\ref{whatdreamsmaycome}.
\item The natural morphism  $\ham\to\mukv$ is smooth and surjective ($\ham$ has a scheme structure from (1), and $\mukv$ from Lemma ~\ref{hamlet}).
\end{enumerate}
\end{lemma}
\begin{proof}
Denote the scheme $\mukv$ by $\mb$. Over this scheme we have a bundle $\mv_{\mathcal{B}}$ and a subbundle $\ms$. Form the total space of  the vector bundle $\mathcal{C}=\shom(\mv_{\mb}/\ms,Q_{\mb})$ over $\mb$ and consider the open subset $U$ of it where the universal morphism $\phi:\mv_{\mc}/\ms_{\mc}\to Q_{\mc}$ is injective (and the cokernel is locally free). It is easy to see that $U\to\mukv$
is smooth and surjective.

Each $\frac{F^j_l(\mv_\mc)}{\ms\cap F^j_l(\mv_{\mc})}$ is a subbundle of $\mv_{\mc}/\ms_{\mc}$ and therefore over $U$,
each $\phi(F^j_a(\mv_{U}))$ is a subbundle of $\mq_{U}$. Hence $H_{\mi}(V,Q,\mk)$ can be represented as the topmost element in a tower of Grassmann bundles over $U$ by successive application of Lemma ~\ref{basic1} (we have to ``choose the flags on $Q$''). This shows that the representing scheme is irreducible and smooth, with a smooth surjective morphism to $\mukv$. It is easy to see that the representing scheme coincides with  $\ham$ at the level of points.
\end{proof}
\bibliographystyle{plain}
\def\noopsort#1{}

\end{document}